\documentclass[12p]{article}
\usepackage{amsmath, amsthm}
\usepackage{amsfonts} %,showkeys}
\usepackage{amssymb}
\usepackage{graphicx}
\usepackage{color}
\usepackage{hyperref}

\newtheorem{thm}{Theorem}[section]

\newtheorem{lem}[thm]{Lemma}
\newtheorem{prop}[thm]{Proposition}
\newtheorem{defn}[thm]{Definition}

\numberwithin{equation}{section}

\newcommand{\RR}{\mathbb{R}}
\newcommand{\ren}{\mathbb{R}^N}
\newcommand{\mS}{{\mathbb S}}
\newcommand{\dx}{\,{\rm d}x}
\newcommand{\dy}{\,{\rm d}y}

\newcommand{\dt}{\,{\rm d}t}
\newcommand{\rd}{{\rm d}}

\def\quotient#1#2{%
    \raise1ex\hbox{$#1$}\Big/\lower1ex\hbox{$#2$}%
}

\def\LL{\mathrm{L}} %per gli spazi L^p
\def\supp{\mathrm{supp}} %per il supporto

\newcommand{\w}{{\Phi}}

\def\ee{\mathrm{e}} %per l'esponenziale
\def\qed{\,\unskip\kern 6pt \penalty 500
\raise -2pt\hbox{\vrule \vbox to8pt{\hrule width 6pt
\vfill\hrule}\vrule}\par}
\definecolor{darkblue}{rgb}{0.05, .05, .65}
\definecolor{darkgreen}{rgb}{0.05, .70, .05}
\definecolor{darkred}{rgb}{0.8,0,0}
\begin{document}
\title{ \bf Optimal Existence and Uniqueness Theory \\ for the Fractional Heat Equation }
\author{Matteo Bonforte$^{\,a}$, %\footnote{e-mail address:~bonforte@calvino.polito.it},
Yannick Sire$^{\,b}$,
~and~ Juan Luis V\'azquez$^{\,a}$} %\footnote{e-mail address:~juanluis.vazquez@uam.es}}
\date{} %%  this cancels date in article format

\maketitle

\begin{abstract}
We construct a theory of existence, uniqueness and regularity of solutions for the fractional heat equation $\partial_t u +(-\Delta)^s u=0$, $0<s<1$,  posed in the whole space $\mathbb{R}^N$ with data in a class of locally bounded Radon measures that are allowed to grow at infinity with an optimal growth rate. We consider a class of nonnegative weak solutions and prove that there is an equivalence between nonnegative data and solutions, which is given in one direction by the representation formula, in the other one by the initial trace. We review many of the typical properties of the solutions, in particular we prove optimal pointwise estimates and new Harnack inequalities.
\end{abstract}

\vspace{1cm}

\small
\parskip 2pt
\linespread{0.7}
\tableofcontents

\normalsize
\parskip 4pt
\linespread{1.0}
\vfill

\hrule
\noindent (a) Departamento de Matem\'{a}ticas, Universidad Aut\'{o}noma de Madrid, Spain\\[2mm]
\noindent (b) Johns Hopkins University, Baltimore, USA, and Universit\'e Aix-Marseille, France

\newpage

\normalsize
\section{Introduction}

We construct a theory of existence, uniqueness, initial traces, as well as a priori estimates and regularity of solutions for the fractional heat equation (FHE) posed in the whole space $\ren$
\begin{equation}\label{FHE}
\partial_t u + (-\Delta)^s u=0, \quad 0<s<1\,.
\end{equation}
The fractional Laplace operator $(-\Delta)^s$ may be defined through its Fourier transform, or by its representation
$$
(-\Delta)^s f(x)=c(N,s)\int_{\ren} \frac{f(x)-f(y)}{|x-y|^{N+2s}}\,\dy
$$
with $0<s<1$, see its basic properties in \cite{Landkoff, Stein}. In the limit $s\to 1$ the standard Laplace operator, $-\Delta$, is recovered (see Section 4 of \cite{DNPV2012}),   but there is a big difference between the local operator $-\Delta$ that appears in the classical heat equation and represents Brownian motion, and the nonlocal family $(-\Delta)^s $, $0<s<1$. These operators are generators of  L\'evy processes that include jumps and long-distance interactions, resulting in anomalous diffusion \cite{App, Bert}.

The main restriction in our paper is that solutions are nonnegative, but they may be quite general otherwise, thus, they are not supposed to be bounded or integrable at any time. We take initial data
\begin{equation}\label{ID}
u(0,\cdot)=\mu_0\,,
\end{equation}
where $\mu_0$ is a nonnegative and locally bounded Radon measure with an admissible growth at infinity
\begin{equation}
\int_{\ren} (1+|x|)^{-(N+2s)}\,d\mu(x)<\infty\,.
\end{equation}
We will prove that such a class, that we call $\mathcal{M}^+_s$, is optimal for existence of solutions of \eqref{FHE}-\eqref{ID}, and it also allows for uniqueness and regularity of the class $\mathcal{W}$ of very weak solutions that we will define in Section \ref{sec.2}. We show that such solutions can be represented by convolution with the so-called fundamental solution in the form
\begin{equation}\label{repres.form.intro}
u(x,t)=\int_{\mathbb R^N} P^t(x-y)\,d\mu_0(y)\,.
\end{equation}
This class of formulas use the fractional heat kernels $P^t(x)=P_s(t,x)$ whose existence is well known. The solutions thus obtained are shown to be very weak and also smooth for all $t>0$. A main issue in that generality is uniqueness. We prove it in Section \ref{sec.uniq} as the result of a delicate Holmgren's argument.

On the other hand, the existence of a unique initial measure for every nonnegative weak solution follows from work of two of the authors \cite{BV2012} in the nonlinear case that we explain in detail here for the present linear equation. In this way we complete the one-to-one correspondence between optimal initial data and very weak solutions in the nonnegative setting.

The question of describing the set of solutions in optimal classes of data is a classical one in the study of the heat equation since the work of Widder \cite{W1, W2} where nonnegative solutions are considered and the optimal class of data is identified as the locally bounded measures that are allowed to grow at infinity in a square quadratic way. The representation formula gives then a unique nonnegative solution that is smooth for positive times and exists globally or locally in time depending on the growth of the initial measure, subcritical or critical respectively.

There are some important differences with the classical heat equation
In that case there was a maximal rate of growth, roughly $u_0(x)\sim e^{a|x|^2}$ for which existence is local in time, up to a blow-up time $T=1/(4a)$ when the solution goes to infinity everywhere in space. No such finite time blow-up seems to happen here according to the established theory. Indeed, the  possible solutions for data with critical growth $u_0(x)\sim |x|^{2s}$ blow up in zero time since
$$
\int \frac{u_0(x)dx}{(1+|x|)^{N+2s}}=\infty.
$$

Concerning the FHE, there are a number of precedents, numerous in the probabilistic theory, cf. \cite{BG1960, BJ2007, ChenKum2008, MK2000, Vald}. Regularity is studied by authors like  \cite{CCV, FK} for equations with more general nonlocal operators, see Section \ref{sect.regularity}. A direct precedent to our paper is \cite{Barrios2014}, where the authors consider nonnegative solutions  with admissible growth at infinity, and establish the validity of the representation formula and a uniqueness theorem in the class of weak solutions as part of the so-called Widder theory. However, measures as initial data are not considered in \cite{Barrios2014} and there is no talk about the existence and uniqueness of initial traces. The  choice of the most convenient class of solutions to describe the ensuing theory is an important issue in that paper.  Here, we take as preferred class the  very weak solutions, which differ from their weak solutions in minor respects, and are suitable for our purpose of obtaining a one-to-one correspondence between solutions and data.

The question of solutions with measure data was considered in the 1980's for a much studied  equation, Porous Medium Equation, a model equation for nonlinear diffusion. Actually, a Widder Theory for that equation was developed in works by B\'enilan et al. \cite{BCP84} and Aronson-Caffarelli \cite{AC83}. The first identifies the class of growing measures for which there is a theory of existence and uniqueness, the second supplies the estimates that allow to identify the class of initial measures of any nonnegative weak solution. Recently, one of the authors \cite{VazBar} constructed the fundamental solution, and Grillo et al. \cite{GMP} addressed the question of uniqueness of solutions with initial data a measure for the fractional porous medium with the restriction that the initial data is a bounded measure. Also, the class of uniqueness is the more restrictive class of weak energy solutions. For some recent work on  fractional nonlinear diffusion equations see \cite{VazNonlin}.

\medskip

\noindent {\sc Organization.} The Widder theory issue contains 3 main theorems that form the basis of the optimality  and equivalence results. We have decided to display the theorems in the following order. All solutions and data are nonnegative unless statement to the contrary. By the label v.w.s. we mean ``very weak solution in the sense of Definition \ref{def.vweak}''.
\begin{itemize}
\item Existence, Thm. \ref{thm.exist}. Existence for data in $\mathcal{M}_s^+$: ``every $\mu_0\in\mathcal{M}_s^+$ generates a v.w.s., which is given by the representation formula, and it takes $\mu_0$ as initial trace in the sense of formula \eqref{init.trace.def}''
\item Initial traces, Thm. \ref{thm.init.trace}. ``Every v.w.s. has a unique initial trace $\mu_0\in\mathcal{M}_s^+$''
\item Uniqueness, Thm. \ref{thm.uniq}. ``Every v.w.s. is uniquely determined by its initial trace''\\ (hence it is given by the formula and corresponds to an initial datum $\mu_0\in\mathcal{M}_s^+$). This completes the equivalence between initial data and solutions given by the solution operator.
\end{itemize}

Note that once solutions exist and are unique, then they are given by the representation formula \eqref{repres.form.intro}, hence comparison is trivially true since the kernel $P^t$ is positive.

After the optimality question is settled, we review a number of topics of the theory following the outline of what is known in the theory of the classical heat equation. Thus, we present a number of results about existence of solutions with $L^1$ data, as well as a number of properties of the solutions, like smoothing effects (also known as ultracontractive estimates). We want to point out the optimal boundedness estimates for general solutions that we develop in Section \ref{sec.boundedness} and the construction of self-similar solutions of \ref{sss}.

We devote some time to comment on the theory for very weak solutions with two signs in Section \ref{sec.signed}.
We point out that in that context solutions may have higher growth rates (both in $x$ and $t$).
The paper ends by a section on comments and open problems and an Appendix.

\medskip

\noindent {\sc Notations.} For all nonnegative functions $f$ and $g$  we use the notation $f \asymp g$
if  there are positive constants $c_1$  and $c_2$ such that $c_1 g(x)\le f(x)\le c_2 g(x)$ in the common domain of definition of both functions. We also write $a\wedge b=\min\{a,b\}$ and
$a \vee  b=\max\{a,b\}$. The Euclidean distance between $x$ and $y$ is denoted by $|x- y|$. We
will use $B_r(x)$ for the open ball centered at $x\in  \ren$ with radius $r > 0.$

\newpage

%%%%%%%%%%%%%%%%%%%%%%%%%%%%%%%%%%%%%%%%%%%%%%%%%%%%%%%%%%%%%%%%%%%%%%%%%%%%%
\section{Some preliminaries and definitions}\label{sec.2}

$\bullet$ The fractional heat equation kernels for $0<s<1$ are  called $P^t(x)$ and are $C^\infty$ functions for all $t>0$ and they are self-similar of the form
$$
P^t(x)=t^{-N/2s}F(|x|t^{-1/2s})
$$
See \cite{BG1960}, \cite{ChenKum2008}. We have
$$
F(\xi) \sim \frac{C}{(1+\xi^2)^{(N+2s)/2}}
$$
Since the Fourier symbol $ e^{-t|\xi|^{2s} }$ is a tempered distribution, it immediately follows  that $P^t \in C^\infty ((0,\infty)\times\ren )$, see  \cite{BJ2007, Kol2000}. Sharp estimates on the behavior of the Heat kernel $P_t$ have been obtained by many authors, see for instance Blumental and Getoor \cite{BG1960}, Chen and Kumagai \cite{ChenKum2008}. They show that
\begin{equation}\label{HK-estimates}
P^t(x)\asymp \frac{1}{t^{N/2s}}\wedge\frac{t}{|x|^{N+2s}}\,.
\end{equation}
A convenient form of the bound would also be
\begin{equation}\label{HK-estimates.1}
P^t(x)\asymp  \frac{t}{\big(t^{1/s}+|x|^2\big)^{(N+2s)/2}}\,.
\end{equation}
The kernel is explicit for $s=1/2$
$$
P_{1/2}(t,x)=c_N\frac{t}{(t^2+| x|^2)^{(N+1)/2}}
$$
A formula and estimates for $F$  are given in \cite{BG1960}.  We will need an interesting estimate for the time derivative,
$$
t\,|\partial_t P^t|/P^t \le C_{N,s}\,,
$$ 
cf. Proposition 2.1 of \cite{VPQR}, that contains also estimates on the behavior of the modulus of the gradient. The lower bound of this inequality is used in Section 6, the optimal constant is $N/(2s)$ and is proved in Lemma \ref{lem.time.monot}.
The fractional heat kernel has a nice smooth Fourier symbol, $\ee^{-|\xi|^{2s}}$, which dictates the finer regularity properties of the solutions ($C^\infty$, analyticity or Gevrey), we refer to Subsection \ref{sec.gevrey} for further details.

\medskip

\noindent\textsc{Weight function. }We define the weight function $\w$ in special form suitable to simplify the technical computations; and we just remark that the main point is that it is $C^2$\,, nonnegative and that $\w\asymp (1+|x|)^{-(N+2s)}$; we therefore define $\w(x)=1$ if $|x|\le 1$ and
\begin{equation}\label{weight}
\w(x):=\dfrac{1}{\left(1+(|x|^2-1)^4\right)^{(N+2s)/8}}\,, \qquad\mbox{if } |x|\ge 1\,.
\end{equation}
Another remarkable property is that there exists a positive constant $c_1$ such that for all $x\in \RR^N$
\begin{equation}\label{weight2}
|(-\Delta)^s \w(x)|\asymp  \w(x)\qquad\mbox{so that}\qquad \left\|\frac{(-\Delta)^s \w}{\w}\right\|_{\LL^\infty(\RR^N)}\le c_1\,.
\end{equation}
The proof of the latter inequality can be checked by a direct calculation, cf. Lemma 2.1 of \cite{BV2012} or Lemma \ref{Lem.phi} in the Appendix.  This weight plays the role of a kind of quasi eigenfunction for the fractional Laplacian on the whole space.
Finally, we can define the weighted $\LL^1$ space as the space of $\LL^1_{\rm loc}(\RR^N)$ functions with finite norm
\[
\|f\|_{\LL^1_\w}:=\int_{\RR^N}|f(x)| \w(x)\dx\,.
\]
We proceed now with the definitions of the preferred class of solutions.
\begin{defn}{\bf Very Weak Solution in $\LL^1_\w$}. \label{def.vweak}
We say that $u$ is a \textsl{very weak solution }to the FHE \eqref{FHE} if (i) $u\in \LL^1_{loc}(0,T: \LL^1_{\w}(\RR^N))$,  and (ii) it satisfies the equality
\begin{equation}\label{def.vweak.int}
\int_0^{+\infty}\int_{\RR^N} u(t,x)\partial_t\psi(t,x)\dx\dt=\int_0^{+\infty}\int_{\RR^N}u(t,x)(-\Delta)^s\psi(t,x)\dx\dt
\end{equation}
for all non-negative $\psi\in C^{\infty}_c((0,+\infty)\times\RR^N)$\,.
\end{defn}

\noindent  {\bf General type of initial data.} We say that a  very weak solution $u(t,x)$ has  initial data $\mu$ if $u\in \LL^1(0,T: \LL^1_{\w}(\RR^N))$,  and there exists a nonnegative  measure $\mu\in \mathcal{M}^+_s$  that is the initial trace of $u$ in the sense that
\begin{equation}\label{init.trace.def}
\int_{\RR^N}\psi\,\rd\mu=\lim_{t\to 0^+}\int_{\RR^N}u(t,x)\psi(x)\,\dx\,,\qquad\mbox{for all }\psi\in C_0(\RR^N)\,.
\end{equation}
Note that since we are only assuming $u\in \LL^1(0,T: \LL^1_{\w}(\RR^N))$, the limit is taken for a.e. $t$. Of course, the solutions given by the convolution  formula \eqref{repres.form} are $C^\infty_{t,x}$ for $t>0$, and we also have  $u\in C((0,T):L^1_\phi(\ren))$. The authors of \cite{Barrios2014} put this last condition in the definition of  weak solution  but we will not do that.

%%%%%%%%%%%%%%%%%%%%%%%%%%%%%%%%%%%%%%%%%%%%%%%%%%%%%%%%%%%%%%%%%%%%%%%

\section{Existence for optimal data. Representation formula}

 The representation formula is defined for all $t>0$, $x\in\ren$ by
\begin{equation}\label{repres.form}
U(t,x)=\int_{\ren} P^t(x-y)\rd\mu_0(y)\,.
\end{equation}
Thanks to the fact that $P^t\in C^\infty ((0,\infty)\times\ren )$ and that its derivatives both in space and in time have a good decay at infinity, it is not difficult to show that also $U\in C^\infty ((0,\infty)\times\ren )$. More precisely, we refer to Section 2 and 3 of \cite{VPQR}.

The first step consists in proving that this formula gives  a very weak solution in the sense of our previous definition.
\begin{thm}[Existence for data in $\mathcal{M}^+_s$]\label{thm.exist}
There exists a nonnegative  very weak solution to the FHE \eqref{FHE} that is given by the representation formula \eqref{repres.form} with initial data $\mu_0 \in \mathcal{M}^+_s$.
\end{thm}

The case of signed data will be discussed in Section \ref{sec.signed}.

\noindent {\sl Proof.~} The proof is split in to several steps.

\noindent$\bullet~$\textsc{Step 1. }The representation formula \eqref{repres.form} is well defined thanks to the decay property of the heat kernel $P^t$ and the matching growth conditions imposed on the data. Since we know that $P^t \in C^\infty ((0,\infty)\times\ren )$ and all the derivatives have also the same (or better) decay properties, we conclude that the solution $U\in C^\infty ((0,\infty)\times\ren )$.

\medskip

\noindent$\bullet~$\textsc{Step 2. }We now show that the formal solution  defined by the representation formula \eqref{repres.form} is indeed a very weak solution in the sense of Definition \ref{def.vweak}. The proof consists just in plugging the formula in the weak formulation \eqref{def.vweak.int} and by checking that all the quantities are finite and equal.\\
We begin with the left-hand side of \eqref{def.vweak.int}:
\begin{equation}\label{thm.exist.01}\begin{split}
\int_0^{+\infty}\int_{\RR^N} U(t,x)\partial_t\psi(t,x)\dx\dt
&=-\int_0^{+\infty}\int_{\RR^N}\int_{\RR^N} \psi(t,x)\partial_t P^t(x-y) \rd\mu_0(y)\dx\dt\\
&=\int_0^{+\infty}\int_{\RR^N}\left(\int_{\RR^N} \psi(t,x)(-\Delta)^s P^t(x-y) \dx\right)\rd\mu_0(y)\dt\\
&=\int_0^{+\infty}\int_{\RR^N}\left(\int_{\RR^N} P^t(x-y) (-\Delta)^s\psi(t,x)  \dx\right)\rd\mu_0(y)\dt\\
\end{split}
\end{equation}
We justify the first calculation by using absolute integrability and Fubini's theorem.
The second and third line then follow by definition of the Heat Kernel, recall that $\partial_tP^t+(-\Delta)^s P^t=0$ pointwise.\\
As for the right-hand side of \eqref{def.vweak.int}, we first rewrite it as follows
\begin{equation}\label{thm.exist.02}
\int_0^{+\infty}\int_{\RR^N}u(t,x)(-\Delta)^s\psi(t,x)\dx\dt
=\int_0^{+\infty}\int_{\RR^N} \left(\int_{\RR^N}P^t(x-y)(-\Delta)^s\psi(t,x) \dx\right)\rd\mu_0(y)\dt\,.
\end{equation}
In order to get equality among \eqref{thm.exist.01} and \eqref{thm.exist.02}\,, we  have to prove that the latter integrals are absolutely convergent, namely we want to show that
\begin{equation}\label{thm.exist.03}
\left|\int_0^{+\infty}\int_{\RR^N}\left(\int_{\RR^N}P^t(x-y)(-\Delta)^s\psi(t,x) \dx\right)\rd\mu_0(y)\dt\right|
\le \int_0^{+\infty}\left|\int_{\RR^N} \w(y)\rd\mu_0(y)\right|\dt<\infty
\end{equation}
To this end, we observe that indeed it is sufficient to prove that for some $c>0$ we have
\begin{equation}\label{thm.exist.04}
\left|\int_{\RR^N}P^t(x-y)(-\Delta)^s\psi(t,x) \dx\right|\le \frac{c}{(1+|y|)^{N+2s}}\asymp \w(y)\qquad\mbox{for all $t>0$\,,}
\end{equation}
so that, in virtue of the assumption $\mu_0\in \mathcal{M}^+_s$ we get \eqref{thm.exist.03}. By a time rescaling argument, it is enough to show that \eqref{thm.exist.04} holds at $t=1$. First, we recall that
\[
P^1(x-y)\asymp \frac{1}{(1+|x-y|)^{N+2s}}\qquad\mbox{and}\qquad |(-\Delta)\psi(x)|\le\frac{c}{(1+|x|)^{N+2s}}
\]
the latter inequality follows by Lemma 2.1 of \cite{BV2012} (see also Lemma \ref{Lem.phi})\,, so that
\begin{equation}\label{thm.exist.05}\begin{split}
\left|\int_{\RR^N}P^1(x-y)(-\Delta)^s\psi(t,x) \dx\right|
\le c\int_{\RR^N}\frac{1}{(1+|x-y|)^{N+2s}}\frac{1}{(1+|x|)^{N+2s}} \dx
\end{split}
\end{equation}
Next, we show that the above integral is bounded by $c\w(y)$ for all $y\in \RR^N$\,, we need to split it in three regions, namely $A_1=B_{|y|/2}(0)$\,, $A_2=B_{|y|/2}(y)$ and $A_3=\RR^N\setminus \{A_1\cup A_2\}$ and estimate the integrals on each region.\\
We begin with $A_1$. We observe that for all $x\in A_1$ we have $|x-y|\ge |y|/2$ so that
\begin{equation}\label{A1}
\int_{A_1}\frac{1}{(1+|x-y|)^{N+2s}}\frac{1}{(1+|x|)^{N+2s}} \dx
\le \frac{2^{N+2s}}{(1+|y|)^{N+2s}}\int_{A_1}\frac{1}{(1+|x|)^{N+2s}} \dx\le \frac{k_1}{(1+|y|)^{N+2s}}\,.
\end{equation}
the latter inequality follows because $(1+|x|)^{-(N+2s)}\in \LL^1(\RR^N)$.\\
As for $A_2$,  we observe that for all $x\in A_2$\,, we have $|x|\ge |y|/2$\,,
\begin{equation}\label{A2}
\int_{A_2}\frac{1}{(1+|x-y|)^{N+2s}}\frac{1}{(1+|x|)^{N+2s}} \dx
\le \frac{2^{N+2s}}{(1+|y|)^{N+2s}}\int_{B_{|y|/2}(y)}\frac{1}{(1+|x-y|)^{N+2s}} \dx\le \frac{k_2}{(1+|y|)^{N+2s}}\,.
\end{equation}
the latter inequality follows because $(1+|x-y|)^{-(N+2s)}\in \LL^1(\RR^N)$ for all $y\in \RR^N$.\\
Finally, the estimate for $A_3$ is the same as for $A_1$\,, indeed we just use that for all $x\in A_3$ we have $|x-y|\ge |y|/2$ and the rest follows similarly to \eqref{A1}.

Summing up, we have proven \eqref{thm.exist.04} which implies \eqref{thm.exist.03}, hence equality among \eqref{thm.exist.01} and \eqref{thm.exist.02} is established.

\noindent$\bullet~$\textsc{Step 3. } We also have to show that the solution $U$ given by the representation formula takes the measure $\mu_0$ as initial trace, in the sense of formula \eqref{init.trace.def}. To this end it is sufficient to notice that the Heat Kernel $P^t$ is an approximation of the identity as $t\to 0^+$\,, hence it converges weakly in the sense of measures to the Dirac Mass $\delta_0$. Since $U(t,x)=P^t * \mu_0$, the previously mentioned weak convergence of the Heat Kernel implies that $U(t,x)\rightharpoonup \mu_0$ as $t\to 0^+$\,, but this is exactly the convergence explicitly stated in formula \eqref{init.trace.def}. The proof of the Theorem is complete.\qed

These solutions will be called for the moment F-solutions, since they are obtained through the representation formula \eqref{repres.form}.

We have already mentioned that the  F-solution of the FHE is  $C^\infty$ smooth in $Q=(0, T )\times \ren$. The regularity at $t=0$ is a big question that heavily depends on the data and that we will try to clarify in what follows. A very simple case happens when $u_0 \in C (\ren)\cap L^\infty(\ren)$, and then it is quite easy to see that the corresponding $F$-solution is continuous down to $t=0$ and also $u\in C([0,\infty):L^\infty(\ren))$. Details about the situation for data $u_0\in \LL^p$ are given in Section \ref{LPtheory}.

%%%%%%%%%%%%%%%%%%%%%%%%%%%%%%%%%%%%%

\section{Initial traces for solutions in weighted spaces}

Our main interest is investigating what happens for more general solutions, to be more specific we pose the \textit{initial trace problem, }that is, wether or not, starting from a nonnegative measure in the class $\mathcal{M}_s^+$\,, very weak solutions in the sense of Definition \ref{def.vweak} will take as initial trace the measure from which we started\,, weakly in the sense of measures, i.e. formula \eqref{init.trace.def}. The positive answer to this question will be given in the next section.

\begin{thm}[\bf Existence and uniqueness of initial traces in $\mathcal{M}^+_s$, \cite{BV2012}]\label{thm.init.trace}
Let $u$ be a nonnegative very weak solution of the FHE \eqref{FHE} in $(0,T]\times\RR^N$. Assume that $\|u(T)\|_{\LL^1_\w(\RR^N)}<\infty$ where $\w$ is as in \eqref{weight}, decaying at infinity as $|x|^{-(N+2s)}$. Then there exists a unique nonnegative Radon measure $\mu\in \mathcal{M}^+_s$ as initial trace, that is
\begin{equation}\label{thm.init.trace.1}
\int_{\RR^N}\psi\,\rd\mu=\lim_{t\to 0^+}\int_{\RR^N}u(t,x)\psi(x)\,\dx\,,\qquad\mbox{for all }\psi\in C_0(\RR^N)\,.
\end{equation}
Moreover, the initial trace $\mu$ satisfies the bound
\begin{equation}\label{thm.init.trace.2}
\int_{\RR^N}\Phi(x)\rd\mu(x)\dx\le \ee^{c_1 T}\,\|u(T)\|_{\LL^1_{\w}(\RR^N)}\,,
\end{equation}
where $c_1>0$ is the constant in the bound \eqref{weight2}.
\end{thm}
The proof has been first given by two of the authors in \cite{BV2012}, and relies on a lemma, Lemma 7.1 of \cite{BV2012}, that gives general conditions for existence and uniqueness of initial traces for very weak solutions of the nonlinear flows corresponding to equation $u_t+ (-\Delta)^s u^m=0$ with $m>0$, i.e. porous medium ($m>1$) or fast diffusion $0<m<1$ and of course FHE, $m=1$. The lemma is then combined with $\LL^1$ weighted estimates (cf. Lemma \ref{lem.weighted.L1}) and gives existence and uniqueness of the initial trace. We will give next a slightly different and self-contained proof of the Theorem, since we feel that it is important here.

\begin{lem}[Weighted $\LL^1$ estimates]\label{lem.weighted.L1}
Let $u$ be a nonnegative very weak solution of the FHE \eqref{FHE}. Let $\w$ be the weight defined as in \eqref{weight}. Then the following inequality holds for all $t,\tau \ge 0$
\begin{equation}\label{lem.weighted.L1.w}
\ee^{-c_1|\tau-t|}\,\|u(\tau)\|_{\LL^1_{\w}(\RR^N)}\le \|u(t)\|_{\LL^1_{\w}(\RR^N)}\le \ee^{c_1|\tau-t|}\,\|u(\tau)\|_{\LL^1_{\w}(\RR^N)}
\end{equation}
where $c_1$ is the constant in the bound \eqref{weight2}.
\end{lem}
\noindent {\sl Proof.~}We first give a formal proof.
\[
\left|\frac{\rd}{\dt}\int_{\RR^d}u(t)\w \dx \right|
= \left|\int_{\RR^d}u(t)\, (-\Delta)^s\w\dx \right|
\le \left\|\frac{(-\Delta)^s\w}{\w}\right\|_\infty\,\int_{\RR^d} u(t)\, \w\dx\le c_1 \int_{\RR^d} u(t)\, \w\dx\,.
\]
Inequality \eqref{lem.weighted.L1.w} then follows then by integration.

The rigorous proof can be done by approximation: take a sequence of test functions $\psi_n(t,x)=\w_n(x)\eta_k(t)\in C^{\infty}_c((0,+\infty)\times\RR^N)$ such that $\eta_k(t)\to \chi_{[\tau_0,\tau_1]}(t)$ and $\partial_t\eta_k(t)\to \delta(\tau_0)-\delta(\tau_1)$ and $\w_n\in C^{\infty}_c(\RR^N)$ with $\w_n\to w$. Then, by Lemma 2.1 of \cite{BV2012} (see also Lemma \ref{Lem.phi} in Appendix) we know that $|(-\Delta)^s\w_n(x)|\le k_1 \w(x)$ and that for $|x|\ge 1$ we have $|(-\Delta)^s\w_n(x)|\ge k_0 \w(x)$ \,, so that we have
\[
\left\|\frac{(-\Delta)^s\w_n}{\w}\right\|_\infty\le c_1\,.
\]
Then plug this sequence of test function in the weak formulation of the equation \eqref{def.vweak.int}\,:
\begin{equation}\label{def.vweak.int.1}\begin{split}
\left|\int_0^{+\infty} \partial_t\eta_k(t) \int_{\RR^N} u(t,x)\w_n(x)\dx\dt\right|
&=\left|\int_0^{+\infty} \eta_k(t)\int_{\RR^N}u(t,x)(-\Delta)^s\w_n(x)\dx\dt\right| \\
&\le \left\|\frac{(-\Delta)^s\w_n}{\w}\right\|_\infty \left|\int_0^{+\infty} \eta_k(t)\int_{\RR^N}u(t,x)\w(x)\dx\dt\right|\\
&\le c_1 \left|\int_0^{+\infty} \eta_k(t)\int_{\RR^N}u(t,x)\w(x)\dx\dt\right|\\
\end{split}\end{equation}
Taking limits as $k\to \infty$ and then as $n\to \infty$ we obtain
\[
\left|\int_{\RR^N} u(\tau_0,x)\w(x)\dx - \int_{\RR^N} u(\tau_1,x)\w(x)\dx\right|
\le c_1 \int_{\tau_0}^{\tau_1} \int_{\RR^N}u(t,x)\w(x)\dx\dt
\]
which by Gronwall lemma gives the desired inequality.\qed

\noindent{\bf Proof of Theorem \ref{thm.init.trace}.~}The proof is divided into three steps.

\noindent$\bullet~$\textsc{Step 1. }\textit{Weighted estimates I. Existence of initial traces. }
Consider the weight function $\w$ defined in \eqref{weight}, and define for $R\ge 1$\,, $\w_r(x)=\w(x/r)$\,, so that $\w_R\equiv 1$ on $B_R(0)$ and recall that it satisfies the estimate \eqref{weight2} in the form
\begin{equation}\label{weight2R}
\left\|\frac{(-\Delta)^s \w_R}{\w_R}\right\|_{\LL^\infty(\RR^N)}\le \frac{c_1}{R^{2s}}\le c_1\,,\qquad\mbox{since we have taken $R\ge 1$}\,.
\end{equation}
Notice that moreover, $\w_R(x)\le R^{N+2s}\w_1(x)$.

We now prove some $\LL^1$-weighted estimates, namely, for all $0\le t\le T_1\le T$ and all $R\ge 1$
\begin{equation}\label{HE.1}\begin{split}
\int_{\RR^d}u(t,x)\w_R(x)\dx
&\le\ee^{c_1(T-t)}\,\int_{\RR^d}u(T,x)\w_R(x)\dx
\le \ee^{c_1\,T}\,\int_{\RR^d}u(T,x)\w_R(x)\dx
\end{split}
\end{equation}
The formal proof of the above inequality is analogous to the proof of Lemma \ref{lem.weighted.L1} and is as follows.
\[
\left|\frac{\rd}{\dt}\int_{\RR^d}u(t)\w_r \dx \right|
= \left|\int_{\RR^d}u(t)\, (-\Delta)^s\w_r\dx \right|
\le \left\|\frac{(-\Delta)^s\w_r}{\w_r}\right\|_{\LL^\infty(\RR^N)}\,\int_{\RR^d} u(t)\, \w_r\dx\le c_1 \int_{\RR^d} u(t)\, \w_r\dx\,.
\]
Then \eqref{HE.1} follows by integration. The rigorous proof can be done by approximation as in the proof of Lemma \ref{lem.weighted.L1}\,, hence we skip the details. As a consequence, we obtain for all $R>0$,
\begin{equation}\label{HE.2}
\int_{B_R(0)}u(t,x)\dx\le \int_{\RR^d}u(t,x)\w_R(x)\dx \le \ee^{c_1\,T}\,R^{N+2s}\int_{\RR^d}u(T,x)\w(x)\dx
:=K_1
\end{equation}
since $\w_R=1$ on $B_R(0)$\,; notice that for $0<R\le 1$ it follows simply by \eqref{HE.1} with the weight $\w=\w_1$\,.

Next, by translation invariance,  it is easy to show that inequality \eqref{HE.2} holds on any ball $B_R(x_0)$\,. Therefore,  we have obtained that for all $R>0$ and all $x_0\in \RR^N$
\begin{equation}\label{HE.2b}
\sup_{t\in (0,T]}\int_{B_R(x_0)}u(t,x)\dx \le K_1  \qquad\mbox{and also}\qquad \limsup_{t\to 0^+}\int_{B_R(x_0)}u(t,x)\dx \le K_1\,.
\end{equation}
this estimate implies weak compactness for measures (to be more precise, weak$^*$ compactness), so that there exists a sequence $t_k\to 0^+$ as $k\to \infty$ with $0<t_k<T_1$\,, and a nonnegative Radon measure $\mu$ so that
\[
\lim_{k\to\infty}\int_{\RR^N}u(t_k,x)\varphi(x)\,\dx=\int_{\RR^N}\varphi\,\rd\mu\qquad\mbox{for all }\varphi\in C^0_c(\RR^N)\,.
\]
The bound on the initial trace: $\mu(B_R(x_0))\le K_1= \ee^{c_1\,T}\,R^{N+2s}\|u(T)\|_{\LL^1_\w(\RR^N)}$ follows from the above bound on the $\limsup$\,.

\noindent$\bullet~$\textsc{Step 2. }\textit{Pseudo-local estimates. Uniqueness. }In order to prove uniqueness of the initial trace we need first to prove the following weighted estimates:
\begin{equation}\label{HE.3}
\int_{\RR^N}u(t,x)\psi(x)\dx\le \int_{\RR^d}u(t',x)\psi(x)\dx
    +K_2 |t-t'|
\end{equation}
for all $0<t,t'\le T_1\le T$ and for all $\psi\in C_c^\infty(\RR^N)$\,. Let us give a formal proof, a rigorous proof can be obtained then by approximation as in the proof of Lemma \eqref{lem.weighted.L1}. Let $\psi\in C_c^\infty(\RR^d)$\,, then by Lemma 2.1 of \cite{BV2012} (see also Lemma \ref{Lem.phi} in Appendix) we know that $|(-\Delta)^s\psi(x)|\le k_1 \w(x)$ and that for $|x|\ge 1$ we have $|(-\Delta)^s\psi(x)|\ge k_0 \w(x)$ \,, so that we have
\[
\left\|\frac{(-\Delta)^s\psi}{\w}\right\|_{\LL^\infty(\RR^N)}\le c_1\,.
\]
Next we calculate
\[
\begin{split}
\left|\frac{\rd}{\dt}\int_{\RR^d}u(t)\psi \dx \right|
&= \left|\int_{\RR^d}u(t)\, (-\Delta)^s\psi\dx \right|
\le \left\|\frac{(-\Delta)^s\psi}{\w}\right\|_{\LL^\infty(\RR^N)}\,\int_{\RR^d} u(t)\, \w\dx
\le c_1 \int_{\RR^d} u(t)\, \w\dx \\
&\le_{(a)} c_1 \ee^{c_1|T-t|}\,\|u(T)\|_{\LL^1_{\w}(\RR^N)}
\le c_1 \ee^{c_1 T}\,\|u(T)\|_{\LL^1_{\w}(\RR^N)}
:= K_2\,.
\end{split}
\]
Notice that in $(a)$ we have used inequality \eqref{lem.weighted.L1.w} of Lemma \eqref{lem.weighted.L1}. Integrating the above differential inequality, we obtain \eqref{HE.3}.

The initial trace whose existence we have proven in Step 1, may depend on the sequence $t_k$\,, hence may not be unique. We will now show  that this is not the case, thanks to estimate \eqref{HE.3}. Assume that there exist two sequences $t_k\to 0^+$ and $t'_k\to 0^+$ as $k\to \infty$\,, so that $u(t_k)\to \mu$ and $u(t'_k)\to \nu$, with $\mu, \nu\in {\cal M}^+(\ren)$. We will prove that
\begin{equation}\label{step.2.uniq.1}
\int_{\RR^N}\varphi\,\rd\mu=\int_{\RR^N}\varphi\,\rd\nu\qquad\mbox{for all }\varphi\in C^{\infty}_c(\RR^N)\,.
\end{equation}
so that $\mu=\nu$ as positive linear functionals on $C^{\infty}_c(\RR^N)$. Then by the Riesz Representation Theorem we know that $\mu=\nu$ also as Radon measures on $\RR^N$. Therefore, it is sufficient to prove \eqref{step.2.uniq.1}\,: estimate \eqref{HE.3} implies that for any $t,t' >0$\,, with $0<t+t'\le T_1\le T$, and any $\varphi\in C^{\infty}_c(\RR^N)$\, we have $|(t+t')-t|=t'$ and
\begin{equation}\label{lem61.step2.1}
\int_{\RR^N}u(t,x)\varphi(x)\dx\le \int_{\RR^N}u(t+t',x)\varphi(x)\dx
    +K_2\,t'\,.
\end{equation}
First we let $t=t_k$ and $t'>0$ to be chosen later, then we let $t_k\to 0^+$ so that $u(t_k)\rightharpoonup \mu$, and we get
\begin{equation}\label{lem61.step2.2}
\int_{\RR^N}\varphi\,\rd\mu\le\int_{\RR^N}u(t',x)\dx
    +K_2\,t'\,.
\end{equation}
Then we put $t'=t'_k$ and let $t'_k\to 0^+$ so that $u(t'_k)\rightharpoonup \nu$ and we obtain the first inequality
\begin{equation}\label{lem61.step2.3}
\int_{\RR^N}\varphi\,\rd\mu\le\int_{\RR^N}\varphi\,\rd\nu\,.
\end{equation}
Then, we proceed exactly in the same way, but we exchange the roles of $t_k$ and $t'_k$ to obtain the opposite inequality $\int_{\RR^N}\varphi\,\rd\mu\ge\int_{\RR^N}\varphi\,\rd\nu\,.$ Therefore we conclude that $\mu=\nu$ as positive linear functionals on $C^{\infty}_c(\RR^N)$\, as desired.

\noindent$\bullet~$\textsc{Step 3. } We still have to pass from test functions $\psi\in C_c^\infty(\RR^d)$ to $\psi\in C_c^0(\RR^d)$  in formula \eqref{thm.init.trace.1}, but this is easy by approximation (mollification). The bound \eqref{thm.init.trace.2} for the initial measure then follows by Lemma \ref{lem.weighted.L1}, namely we get
\[
\int_{\RR^N}\Phi(x)\rd\mu(x)\dx\le \ee^{c_1 T}\,\|u(T)\|_{\LL^1_{\w}(\RR^N)}\,.
\]
The proof is concluded.\qed

 %%%%%%%%%%%%%%%%%%%%%%%%%%%%%%%%%%%%%%%%%%%%%%%%%%%%%%%%%%%%%%%%%%%%%%%%%%%%
\section{Uniqueness of very weak solutions for optimal data}\label{sec.uniq}

\begin{thm}\label{thm.uniq} Every nonnegative very weak solution in the sense of Definition $\ref{def.vweak}$ is uniquely determined by its initial trace $\mu_0\in \mathcal{M}_s^+$.
\end{thm}

\noindent {\sl Proof.~}We split the proof into several steps.

\noindent$\bullet~$\textsc{Step 1. }Fix $T>0$. We recall the weak formulation of our problem, according to Definition \ref{def.vweak}, namely
\begin{equation}\label{uniq.step1.1}
\int_0^{T}\int_{\RR^N} u(t,x)\partial_t\psi(t,x)\dx\dt
=\int_0^{T}\int_{\RR^N}u(t,x)(-\Delta)^s\psi(t,x)\dx\dt
\end{equation}
for every  test function $\psi\in C^{\infty}_c((0,T)\times\RR^N)$. The main idea towards uniqueness is to show that the difference of two nonnegative solutions $u_1-u_2$ corresponding to the same nonnegative initial trace\,, is indeed the zero solution. We have now to be careful that letting $u=u_1-u_2$\,, $u$ is not anymore a nonnegative solution, but it is clear that $u_+=u_1\ge 0$ and $u_-=u_2\ge 0$\,, so that $|u(t,\cdot)|=u_1(t,\cdot)+u_2(t,\cdot)\in \LL^1_\w(\RR^N)$ for all $t>0$.

To prove uniqueness, i.e. to prove that indeed $u=0$ for all $t>0$\,, we would like to use as test function the solution $\varphi(t,x)$ of the dual problem
\begin{equation}\label{backward.HE}
-\partial_t\varphi(t,x)+(-\Delta)^s\varphi(t,x)=f(t,x)\in C^{\infty}_c((0,T)\times\RR^N)
\end{equation}
 which is a backward fractional  heat equation: We take final data
\begin{equation}\label{backward.HE.data}
\varphi(T,x)= 0 \,,
\end{equation}
and use Duhamel's formula to represent the solution
$$
\varphi(t)=P^{T-t}\,\varphi(T)+\int_t^T P^{\tau-t}\ast f(\tau)\,d\tau\,.
$$
where $0<t<T$ and we use the short  notation $\varphi(t)=\varphi(t,\cdot)$. In view of the zero final data this simplifies to
$$
\varphi(t,\cdot)= \int_t^T P^{\tau-t}\ast f(\tau,\cdot)\,\rd\tau\,.
$$
Since $f$ is smooth and compactly supported, we can conclude that $\varphi(t,x)$ is smooth and decreases at infinity like $(1+|x|)^{-(N+2s)}$ uniformly in $0<t<T$. Indeed, we recall that
\begin{equation}\label{HK-estimates.1b}
P^t(x)\asymp  \frac{t}{\big(t^{1/s}+|x|^2\big)^{(N+2s)/2}}\,.
\end{equation}
and by a scaling argument (and the fact that $f$ is compactly supported in time)\,, it is sufficient to show that for all $|x|$ large enough we have that
\begin{equation}\label{uniq.1.1}
\left[P^{1}\ast f(1,\cdot)\right](x)\asymp \frac{1}{\big(1+|x|^2\big)^{(N+2s)/2}}
\end{equation}
Without loss of generality we can assume that $\supp(f(1,\cdot))\subset B_{R_0}(0)=B_{R_0}$ for some $R_0>0$\,, and we will write $f(x)$ instead of $f(1,x)$ with a little abuse of notation. We then consider $|x|\ge 2 R_0\ge 2 |y|$. We begin with the lower bounds and we observe that for all $y\in B_{R_0}$ and all $|x|\ge 2 R_0$ we have $|x-y|\le |x|+|y|\le |x|+R_0\le 2|x|$ so that
\begin{equation}\label{uniq.1.2}\begin{split}
\left[P^{1}\ast f\right](x)&\ge c_0\int_{B_{R_0}}\frac{f(y)}{\big(1+|x-y|^2\big)^{(N+2s)/2}}\dy\\
&\ge \frac{c_0}{\big(1+(2|x|)^2\big)^{(N+2s)/2}}\int_{B_{R_0}}f(y)\dy \ge \frac{k_0}{\big(1+|x|^2\big)^{(N+2s)/2}}\,.
\end{split}
\end{equation}
Now we prove the upper bounds, noticing that since $|y|\le R_0$ and $|x|\ge 2 R_0\ge 2|y|$ we have $|x-y|\ge |x|-|y|\ge |x|/2$ so that
\begin{equation}\begin{split}
\left[P^{1}\ast f\right](x)&\le c_1\int_{B_{R_0}}\frac{f(y)}{\big(1+|x-y|^2\big)^{(N+2s)/2}}\dy\\
&\le \frac{c_1}{\big(1+(|x|/2)^2\big)^{(N+2s)/2}}\int_{B_{R_0}}f(y)\dy \le \frac{k_1}{\big(1+|x|^2\big)^{(N+2s)/2}}\,.
\end{split}
\end{equation}
By rescaling we obtain \eqref{HK-estimates.1b}.

\noindent$\bullet~$\textsc{Step 2. }Taking as test function $\psi=\varphi$ in formula \eqref{uniq.step1.1} is not allowed since $\varphi$ is not compactly supported, neither in space nor in time near $t=0$. Therefore, we take $\psi(t,x)=\varphi(t,x)\zeta_R(x)\theta_n(t)$ where $\zeta_R$ does a cutoff in the space variable and $\theta$ is a cutoff in time near $t=0$. More precisely we choose $\theta_1$ a smooth function that grows from $\theta_1(0)=0$ to $\theta_1(1)=1$ and we define $\theta_n(t)=\theta_1(nt)$. As for $\zeta_R(x)$\,, we consider the standard cutoff function which takes value 1 on the ball $B_R$ and vanishes outside $B_{2R}$.
Then we have
$\partial_t\psi=\zeta_R\theta_n\partial_t\varphi + \varphi\, \theta'_n\zeta_R$ and by the well-known product formula, we get
$$
(-\Delta)^s\psi=\theta_n \,\left[\zeta_R(-\Delta)^s\varphi+ \varphi(-\Delta)^s\zeta_R+ B(\varphi, \zeta_R )\right]\,,
$$
where the bilinear form $B$ is defined as follows
\[
B(\varphi, \zeta_R )(x)=\int_{\RR^N}\frac{(\varphi(x)-\varphi(y))(\zeta_R(x)-\zeta_R(y))}{|x-y|^{N+2s}}\dy\,.
\]
We then plug $\psi$ into the weak formulation \eqref{uniq.step1.1} and we use the fact that $\varphi$ solves the backward heat equation to get
\begin{equation}\label{uniq.2.1}
\int_0^{+\infty}\int_{\RR^N} u(t,x)f(t,x)\zeta_R(x)\theta_n(t)\dx\dt
+I_1+I_2= \int_0^{+\infty}\int_{\RR^N}u(t,x)\varphi(t,x)\zeta_R(x)\theta'_n(t)\dx\dt
\end{equation}
with
$$
I_1:=\int_0^{+\infty}\int_{\RR^N} u(t,x)\varphi(t,x)(-\Delta)^s\zeta_R(x)\theta_n(t)\,\dx\dt\,,
$$
and
$$
I_2:=\int_0^{+\infty}\int_{\RR^N} u(t,x)B(\varphi(t,x),\zeta_R(x))\theta_n(t)\,\dx\dt\,.
$$
We also define
\[
J:=\int_0^{+\infty}\int_{\RR^N}u(t,x)\varphi(t,x)\zeta_R(x)\theta_n'(t)\dx\dt\,.
\]
The idea is to prove that in the limit $R\to\infty$, $n\to\infty$ the three terms $I_1$\,, $I_2$ and $J$ go to zero, (since $\zeta_R\to 1$, $\theta_n(t)\to 1$) so that formula \eqref{uniq.2.1} becomes
$$
\int_0^{+\infty}\int_{\RR^N} u(t,x)f(t,x)\dx\dt=0
$$
and then, by the test Lemma we would conclude that $u\equiv 0$, the above expression is zero for all test functions $f\in C^{\infty}_c((0,+\infty)\times\RR^N)$\,. This would conclude the proof of uniqueness.

Therefore, it only remains to prove that in the limit $R\to\infty$, $n\to\infty$ the three terms $I_1$\,, $I_2$ and $J$ go to zero, and this will be done in the next Step.

\noindent$\bullet~$\textsc{Step 3. }The important part is the integral J, where we have to use the form in which solutions approach the initial measures. As far as $I_1$ and $I_2$ are concerned, up to the smooth cut-off function $\theta_n$, these integrals have been considered in \cite{Barrios2014}, Section 2 ($I_1$ formula (2.14) and $I_2$ formula (2.33)), the proof of convergence is the same in our case and we do not reproduce it here.

\noindent To prove that $J\to 0$ as $n\to\infty$\,, we split the integral $J$ in two parts, $J=J_{1,n}+J_{2,n}$ where
$$
J_{1,n}=\int_0^{1/n}\theta_n'(t)\int_{\ren}u(t,x)\varphi(0,x)\zeta_R(x)\,\dx\dt
$$
and
$$
J_{2,n}=\int_0^{1/n}\theta_n'(t)\int_{\ren}u(t,x)\varphi(0,x)
\frac{\varphi(t,x)-\varphi(0,x)}{\varphi(0,x)}\zeta_R(x)\,\dx\dt
$$
We begin with $J_{1,n}$. Let us first recall that the initial trace of $u=u_1-u_2$ is $\mu=0$\,, because by definition of initial trace we have for $i=1,2$ that
\begin{equation}\label{uniq.3.1}
\lim_{t\to 0^+}\int_{\RR^N}u_i(t,x)\phi(x)\,\dx=\int_{\RR^N}\phi\,\rd\mu_0\,,\qquad\mbox{for all }\phi\in C_0(\RR^N)\,,
\end{equation}
because both $u_1$ and $u_2$ have the same initial trace $\mu_0$. Hence
\begin{equation}\label{uniq.3.2}
\lim_{t\to 0^+}\int_{\RR^N}u(t,x)\phi(x)\,\dx=0\,,\qquad\mbox{for all }\phi\in C_0(\RR^N)\,.
\end{equation}
The above formula implies that the function $t\mapsto Y(t):=\int_{\RR^N}u(t,x)\phi(x)\,\dx$ is a continuous function on $(0, T/2)$ and it is right-continuous at $0$ and $Y(0^+)=0$. Moreover, we have that $\theta_n'(t)\chi_{[0,1/n](t)}\to \delta_0(t)$ in the weak topology of measures. As a consequence\,,
\[
\lim_{n\to\infty}J_{1,n} =\lim_{n\to\infty}\int_0^\infty \theta_n'(t)\chi_{[0,1/n](t)} Y(t)\dt=Y(0^+)=0\,,
\]
which proves that $J_{1,n}\to 0$ uniformly in $R$.

Next we deal with $J_{2,n}$. We first observe that for all $t\in [0,T/2]$ and all $x\in \RR^N$ we have
\begin{equation}\label{uniq.3.2b}
\left|\frac{\varphi(t,x)-\varphi(0,x)}{\varphi(0,x)}\right|\le Ct\,.
\end{equation}
We only have to check the above inequality for large values of $|x|$. We have proven in Step 1 that for large $|x|$, we have $\varphi(t,x)\asymp c(t)\big(1+|x|^2\big)^{-(N+2s)/2}$\,, for some bounded function $c(t)$; in particular $\varphi(t,x)$ and $\varphi(0,x)$ has the same behavior for large $|x|$ \,, so that we can conclude that
\[
\frac{\varphi(t,x)}{\varphi(0,x)}\le k_3\, c(t)\le k_4\qquad\mbox{for all $t\in [0, T/2)$ and all $x\in \RR^N$\,.}
\]
Moreover, $\varphi(t,x)$ is a smooth function with respect to the $t$ variable\,, so that by the mean value Theorem
\[
|\varphi(t,x)-\varphi(0,x)| \le \varphi(\tilde{t},x)t\qquad\mbox{for all $t\in [0, T/2)$ and suitable $\tilde{t}\in (0,t)$.}
\]
Therefore, combining the two above estimates we get for all $x\in \RR^N$
\[
\left|\frac{\varphi(t,x)-\varphi(0,x)}{\varphi(0,x)}\right|\le \frac{\varphi(\tilde{t},x)}{\varphi(0,x)} t\le k_4 t
\]
that is \eqref{uniq.3.2b}.
Next, we prove that
\begin{equation}\label{uniq.3.3}
\int_{\ren}|u(t,x)|\varphi(0,x)\zeta_R(x)\,\dx\le k_5<\infty\qquad\mbox{for all $t\in [0,T]$}\,.
\end{equation}
Indeed, we know by Lemma \ref{lem.weighted.L1} that any weak solution $u_i$ satisfies the estimate
\begin{equation}\label{lem.weighted.L1.w.b}
\|u_i(t)\|_{\LL^1_{\w}(\RR^N)}\le \ee^{c_1\,T}\,\|u_i(T)\|_{\LL^1_{\w}(\RR^N)}\,,\qquad\mbox{for all }t\in[0,T]\,.
\end{equation}
and we also know that $\|u_i(T)\|_{\LL^1_{\w}(\RR^N)}\le k_6<\infty$\,, so that, recalling that $\varphi(0,x)\asymp\w$
\[\begin{split}
\int_{\ren}|u(t,x)|\varphi(0,x)\zeta_R(x)\,\dx
&\le  \int_{\ren}u_1(t,x)\varphi(0,x) \,\dx +\int_{\ren}u_2(t,x)\varphi(0,x) \,\dx\\
&\le  \int_{\ren}u_1(t,x)\w \,\dx +\int_{\ren}u_2(t,x)\w \,\dx\\
&\le \ee^{c_1\,T}\left(\|u_1(T)\|_{\LL^1_{\w}(\RR^N)}+\|u_2(T)\|_{\LL^1_{\w}(\RR^N)}\right)\le 2k_6\ee^{c_1\,T}=k_7<\infty
\end{split}\]
Finally we can show that $J_{2,n}\to 0$ as $n\to \infty$ uniformly in $R$\,, indeed
\[\begin{split}
|J_{2,n}|&=\left|\int_0^{1/n}\theta_n'(t)\int_{\ren}u(t,x)\varphi(0,x)\frac{\varphi(t,x)-\varphi(0,x)}{\varphi(0,x)}\zeta_R(x)\,\dx\dt\right|\\
&\le \int_0^{1/n}|\theta_n'(t)| \int_{\ren}|u(t,x)|\varphi(0,x)\left|\frac{\varphi(t,x)-\varphi(0,x)}{\varphi(0,x)}\right|\zeta_R(x)\,\dx\dt\\
&\le k_4 \int_0^{1/n}|\theta_n'(t)|t \int_{\ren}|u(t,x)|\varphi(0,x)\zeta_R(x)\,\dx\dt\\
&\le k_4 k_5 \int_0^{1/n} t\theta_n'(t) \dt = k_4 k_5\left[ \frac{\theta_n(1/n)}{n} - \int_0^{1/n}\theta(t)\dt\right]\xrightarrow{n\to \infty}0\\
\end{split}
\]
where we have used \eqref{uniq.3.2b} , \eqref{uniq.3.3} and the fact that $\theta'\ge 0$ since $\theta$ is growing. We have also used that $\theta_n$ is bounded, in the last step.

\noindent The proof of the uniqueness Theorem is concluded.\qed

\medskip

\noindent{\bf Historical note.} This type of uniqueness proof based on duality is very typical of linear partial differential equations with real analytic coefficients with the name of Holmgren's theorem, see \cite{Holm} and the books by Smoller or Treves. Wikipedia says: \sl In the theory of partial differential equations, Holmgren's uniqueness theorem, or simply Holmgren's theorem, named after the Swedish mathematician Erik Albert Holmgren (1873-1943), is a uniqueness result for linear partial differential equations with real analytic coefficients. \rm
As we see, the original proof was rather strong in assumptions.

  The technique has been extended to nonlinear equations by experts in nonlinear PDES and the assumptions have been relaxed. Thus,  it was used by Kamin in
1961 for the Stefan Problem, and then by Kalashnikov in 1963 for the porous medium equation in 1d with growing data; both were Olga Oleinik's students. The proof in all dimensions for the PME is due to B\'enilan, Crandall and Pierre in their famous paper in 1984. It is described in Theorem 6.5 of Vazquez's book PME.

  In this paper we extend to optimal measures as initial data the type of proof done for the
fractional heat equation in the paper by Barrios et al  in 2014 \cite{Barrios2014}.

\

%%%%%%%%%%%%%%%%%%%%%%%%%%%%%%%%%%%%%%%%%%%%%%%%%%%%%%%%%%%%%%%%%%%%%%%

\newpage

\section{$\LL^1$ and $\LL^p$ theory}\label{LPtheory}

This section contains material that is more or less known or used in the literature, and we
present it for completeness, asking the reader to supply missing details.

\noindent$\bullet$ {\bf Contraction properties. } The fractional operator is a maximal monotone operator, subdifferential of a convex functional in $L^2(\ren)$,  hence it generates a contraction semigroup in that functional space. The contractivity property extends to  all $\LL^p(\RR^N)$ spaces, $p\in [1,\infty]$\,, namely
\[
\|u(t)-v(t)\|_{\LL^p(\RR^N)}\le \|u_0-v_0\|_{\LL^p(\RR^N)}
\]
This is a consequence of the fact that $(-\Delta)^s$ is $m$-accretive in all $L^p$ spaces.
The contractivity can also be explained in a simple way, as follows:
\[
\frac{\rd}{\dt}\int_{\RR^N}|u|^p(t,x)\dx=-p\int_{\RR^N}|u|^{p-1}{\rm sign}(u)(-\Delta)^s u\dx \le 0
\]
Indeed, to ensure the positivity of the last term, we need the so-called Stroock-Varopoulos inequality, namely
\[
\int_{\RR^N}u^{p-1}(t,x)(-\Delta)^s u(t,x)\dx \ge \frac{4(p-1)}{p^2} \int_{\RR^N}\left|(-\Delta)^{\frac{s}{2}} u^{\frac{p}{2}}(t,x)\right|^2\dx
\]
see \cite{Var} and also \cite{DPQRV2}. Of course the above proof needs to be justified, by approximation, hence the need for a slightly more general inequality then the one above, cf. for instance Lemma 5.2 of \cite{DPQRV2}. The Stroock-Varopoulos inequality is also the basis to begin a Moser iteration based on fractional Gagliardo-Nirenberg inequalities, that would prove smoothing effects without the use of the heat kernel, a method that allows to deal with more general diffusion operators.

In this way, a semigroup of contractions is obtained in all $L^p$ spaces $p\ge 1$\,, indeed, by writing $S_t(\mu_0)=P^t\ast u_0$, then it is easy to check that
\begin{equation}
 S_t\circ S_{t'}(u_0)=S_{t+t'}(u_0),\qquad\mbox{and that}\qquad S_0(u_0)=u_0,
\end{equation}
since we recall that $P^t$ is an approximation of the identity (as $t\to 0^+$) and that $\lim\limits_{t\to 0^+}P^t=\delta_0$ weakly in the sense of measures\,, as carefully explained in Step 3 of the above Existence Theorem \ref{thm.exist}.

\medskip

\noindent $\bullet$ {\bf  Time monotonicity}

\begin{lem}\label{lem.time.monot}
Let $u$ be a nonnegative very weak solution of the FHE \eqref{FHE}. Then the time derivative $u_t$ satisfies the following inequality:
\begin{equation}\label{time.monotonicity}
u_t(t,x)\ge -\frac{N}{2s}\frac{u(t,x)}{t}\qquad\mbox{for all $t>0$ and all $x\in \RR^N$.}
\end{equation}
that implies that the function $t\mapsto t^{N/2s}u(t,x)$ is monotone nondecreasing for all $x\in \RR^N$.
\end{lem}
\noindent {\sl Proof.~}We first recall that the very weak solution $u(t,x)$ is unique and can be expressed in terms of the representation formula in terms of its initial trace, as consequence of Theorems \ref{thm.exist}, \ref{thm.uniq} and \eqref{thm.init.trace}. Therefore
\begin{equation}\label{repres.form.1}
u(t,x)=\int_{\ren} P^t(x-y)\rd\mu_0(y)\,.
\end{equation}
where $P^t(x)=t^{-N/2s}F(|x|t^{-1/2s})$, where $F:\RR\to\RR$ is a smooth decreasing function. Moreover,
\[
\partial_t P^t(x)= -\frac{N}{2s}\frac{F(|x|t^{-1/2s})}{t^{\frac{N}{2s}+1}}
-\frac{|x|}{2s}\frac{F'(|x|t^{-1/2s})}{t^{\frac{N+1}{2s}+1}}
\ge -\frac{N}{2s\, t}P^t(x)
\]
since $F'\le 0$. As a consequence,
\begin{equation}\label{repres.form.2}
\partial_t u(t,x)=\int_{\ren} \partial_tP^t(x-y)\rd\mu_0(y)\ge  -\frac{N}{2s\, t} \int_{\ren} P^t(x-y)\rd\mu_0(y)= -\frac{N}{2s}\frac{u(t,x)}{t}.
\end{equation}
The monotonicity of  $t\mapsto t^{-N/2s}u(t,x)$  follows by integration. This concludes the proof. \qed

\medskip

\noindent $\bullet$ {\bf $\LL^1-\LL^\infty$ Smoothing Effect. Ultracontractivity}

\begin{thm}
Let $u$ be a very weak solution of the FHE \eqref{FHE} corresponding to the nonnegative initial datum $0\le u_0\in \LL^1(\Omega)$. Then there exist a constant $K_p>0$ depending only on $p\ge 1$\,, $0<s\le 1$ and $N$\,, such that for all $t>0$
\begin{equation}\label{thm.smooth.eff.std}
\|u(t)\|_{\LL^\infty(\RR^N)}\le \frac{K_p}{t^{\frac{N}{2sp}}}\|u_0\|_{\LL^p(\RR^N)}\,.
\end{equation}
\end{thm}
\noindent {\sl Proof.~}The case $p=1$ uses the representation formula and the inequality $P^t(t,x)\le C t^{-\frac{N}{2s}}$. The case $p>1$ follows by interpolation.\qed

\noindent $\bullet$ {\bf Asymptotic behaviorfor $L^1$ data.} The following theorem parallels what happens for the heat equation,

\begin{thm} Let $u(t,x)$ be the very weak solution with initial data $u_0\in L^1(\ren)$, $u_0\ge 0$, and let $M=\int u_0(x)\,dx$ be its mass, which is constant in time. Then as $t\to\infty$ we have
\begin{equation}
\lim_{t\to\infty}t^{\alpha_p}\|u(t,x)-MP^t(x)\|_{L^p(\ren)}=0
\end{equation}
with $\alpha_p=\frac{n(p-1)}{2sp}.$
\end{thm}

\noindent {\sl Proof.} Write the difference as
$$
u(t,x)-MP^t(x)=\int u_0(y)(P^t(x-y)-P^t(x))\,dx
$$
and proceed as in the proof for the heat equation both in the case
of the $L^1$ norm and in the case of the $L^\infty$ norm. The rest of the $L^p$ norms
are obtained by interpolation.  We leave the rest as an exercise, since this is not difficult. \qed

A similar result is also true for the porous medium equation, the $p$-Laplacian equation and other nonlinear diffusion equations, but the proofs need nonlinear techniques, cf. \cite{KamVaz, VazPME}.

\medskip

\noindent$\bullet~$ {\bf Symmetrization}. In the theory of elliptic and parabolic equations, the technique of Schwarz symmetrization is a frequently used  tool to obtain a priori bounds for classical and weak solutions in terms of general information on the data, see \cite{Ba76, Ta79} from the huge literature.

We define the symmetric rearrangement of a set $A\subset \RR^n$ as the ball centered at $x=0$ with the same volume, that we denote by $A^*$. We also define the symmetric decreasing rearrangement of an integrable function $f\in L^1(\RR^n)$, that we denote by $f^*$; briefly stated, $f^*$ is a radially symmetric function defined in $\RR^n$, it is decreasing as a function of $r=|x|$, and its level super-sets $\{x: f^*(x)>\lambda>0\}$ have the same measure as those of $f$.

The application is as follows; In the elliptic theory we want to compare the solution
$u$ of a  problem with right-hand side $f$ with the solution $v$ of the problem with rearranged right-hand data $f^*$. In the parabolic case symmetrization is applied only on the space variable and then we want to compare $u(t,x)$ with $v(x,t)$ for all fixed times $t>0$. The specific theorems may be pointwise comparison  of integral comparison. The net benefit is that solving radial problems is easier, and even in some times there are explicit formulas.

In our present linear setting the following comparison result is easy to derive.

\begin{thm} Let $u(t,x)$ be the very weak solution with initial data $u_0\in L^1(\ren)$, $u_0\ge 0$, and let $v(t,x)$ be the very weak solution with initial data $v_0=u_0^+$. Then for every $t>0$ we have
$$
0\le u(t,x) \le \|v(t)\|_\infty=v(t,0)\,.
$$
\end{thm}

\noindent {\sl Proof.} We use the  well-known Hardy-Littlewood Lemma \cite{HLP}  inequality
$$
\int fg\,dx\le \int f^*g^*\,dx
$$
on the solution given by the convolution formula to get
$$
\int u_0(y)P^t(x-y)\,dy\le \int u_0^*(y) (P^t(x-y))^*\,dx\,.
$$
Now, it is immediate that $P^t(x-y)^*=P^t(y)$ and we get the relation
$u(t,x)\le v(t,0)$, which allows to compare the $L^\infty$ norms of $u(t)$ and $v(t)$.

Stronger results can be obtained by using the full machinery of symmetrization, which also applies to nonlinear fractional flows, cf. \cite{VazVol1, VazVol2}, but the proofs are much longer.

%%%%%%%%%%%%%%%%%%%%%%%%%%%%%%%%%%%%%%%%%%%%%%%%%%%
\section{Quantitative local boundedness and optimal existence}
\label{sec.boundedness}

\subsection{The  weighted estimate}

We want to estimate the behaviorof the constructed solution $u=P^t\ast \mu_0$ for $t>0$ and prove that it is a locally bounded function of $x$ with precise estimates. Here is the first result.

\begin{thm}\label{Thm.upper.gen} Let $u=S_t\mu$ the very weak solution with initial measure $\mu_0\in \mathcal{M}^+_s$ and let $\|\mu_0\|_\Phi:= \int_{\RR^N}\w\rd\mu_0$. There exists a constant $C(N,s)$  such that for every $t>0$ and $x\in\ren$
\begin{equation}\label{quant.est}
u(t,x)\le C\,\|\mu_0\|_\Phi(t^{-{N/2s}}+t) (1+|x|)^{N+2s}, \qquad
\end{equation}
\end{thm}

\noindent {\sl Proof.} As an opening, we examine the simplest  case where $x=0$ and  time is bounded and bounded away from zero. We have
$$
u(t,0)=\int_{\ren}P^t(y)\,d\mu(y)
$$
Under our assumptions we have $P^t(y)\asymp \Phi$, hence $u(t,0)\asymp \|\mu_0\|_\Phi$. The same argument works when $x$ is bounded.

\medskip

\noindent$\bullet$ In other to proceed further, it will be convenient to decompose $\mu_0=\mu_1+\mu_2$ where $\mu_1$ is the part of $\mu_0$ located in the open ball of radius 1. Then we write the solution as $u_1=S_t (\mu_1)$ and $u_2=S_t (\mu_2)$. We know the behaviorof the fist by simple inspection of the kernel: $u_1(t)$ is a bounded function for every $t>0$ and
$$
\|u_1(t)\|_\infty\le  C\,t^{-N/2s}\|\mu_1\|_M,
$$
where $\|\cdot\|_M$ is the standard norm for bounded measures. Moreover, the estimate is sharp for all $t>0$ since the fundamental solution satisfies it.

\medskip

\noindent$\bullet$ Next, we examine the case $t\sim 0$ for $x\in B_{1/2}$. In view of the previous paragraph, we only need to consider $u_2$.
For $|y|>1$ we have
$$
P^t(x-y)\sim C\frac{t}{(|x-y|^2+t^{ 1/s})^{(N+2s)/2}}\le C_1\frac{t}{(|y|^2+1)^{(N+2s)/2}}\,,
$$
so that
$$
u_2(t,x)\le C_1 t \|\mu_2\|_\Phi\,.
$$
Adding the estimates for $u_1$ and $u_2$ we get the estimate
$$
u(t,x)\le C_2  t^{-N/2s} \|\mu\|_\Phi\,,
$$
valid for all $t\le 1$ and $|x|\le 1/2$. A simple scaling gives the proof for $|x|\le R$, and then
$C_2 $ depends also on $R$. This far, we have proved that for all bounded $x$ and $0<t< \infty$ we have
\begin{equation}
u(t,x)\le C\,(t^{-{N/2s}}+t)\|\mu_0\|_\Phi\,,
\end{equation}
where $C(R,N,s)$.

\medskip

\noindent$\bullet$ Next, we examine the case $t\gg 1$ and $x$ bounded. Estimating the expression for $P_t$
like before we get $u_2(t,x), u(t,x)\le C_1 t \|\mu_2\|_\Phi\,$.

\medskip

\noindent {\bf  Comment. }  Separating the contributions for $|y|\le t^{1/2}$ from $|y|\ge t^{1/2}$ in the last paragraph, we may even arrive at $u_2(t,x), u(t,x)=o(t)$ as $t\to\infty$.
\medskip

\noindent$\bullet$ To conclude the proof, we need to translate this into the growth in $|x|$ large for finite fixed $t>0$. We may use scaling
which seems amazing. We change the coordinates $x'=x-a$ with $|a|>1$ and consider the previous formula
in the new coordinates. Then, putting tildes for the functions in the displaced coordinates,
$$
u(t,a)=\widetilde{u}(t,0)\le C\,(t^{-{N/2s}}+t)\|\widetilde{\mu}_0\|_\Phi\,.
$$
Now, the comparison of the weights shows that $\|\widetilde{\mu}_0\|_\Phi\le a^{N+2s}\|\mu_0\|_\Phi$. Hence
$$
u(t,a)\le C\,\|\mu_0\|_\Phi(t^{-{N/2s}}+t) (1+|a|)^{N+2s}\,.
$$
This is the final result. \qed

\medskip

\noindent {\bf Optimality of the exponents. }The exponent in the $x$ variable cannot be improved
in view of the following example. Let $x_0=a e_1$ with $a>0$ large and take as initial data
$$
\mu_a(x)=c\, a^{N+2s}\delta(x-x_0)
$$
then it is easy to see that $\|\mu_a\|_\Phi \asymp 1$ with a constant that does not depend on $a=|x_0|\gg 1$. On the other and we know that the corresponding solution satisfies
$$
u_a(t,x_0)=c\, |x_0|^{N+2s}\,t^{-N/2s}, \qquad
$$
We recall that $|x_0|$ is as large as we wish.

For the optimality of the $t$ exponents we refer to next section.

%%%%%%%%%%%%%%%%%%%%%%%%%%%%%%%%%%%%%%%%%%%
\subsection{Existence of self-similar solutions in classes of growing data}\label{sss}

We consider initial data of the form
$$
U_0(x)=C(\theta)|x|^{\gamma}
$$
where $0\le \gamma<2s$ and $C(\theta) \ge 0$ is a continuous function of $\theta\in \mS^{N-1}$. Such a function
falls into the admissible initial data $\mathcal{M}_s^+(\ren)$. It follows from our results that there
is a unique very weak solution $U(x,t)$ and we have a number of properties:

(i) $U$ is positive and $C^\infty$ in $Q=(0,\infty)\times \ren$\,,

(ii)  We have a priori estimates from the previous section.

(iii) If $C(\theta)$ is constant, then $U$ is radially symmetric with respected to the space variable, and  increasing in $|x|$ for fixed $t>0$,

Next, we take into account the invariance of the initial data under the scaling $T_\lambda (U_0)=\lambda^{-\gamma} U_0(\lambda x)$. Using the equation we conclude that
$$
U_\lambda(t,x)=\lambda^{-\gamma} U(\lambda^{2s } t, \lambda x)
$$
is another solution with the same initial data, hence  $U_\lambda(t,x)=U(t,x)$
which means that
$$
U(t,x)=t^{\gamma/2s}  F(x\,t^{-1/2s})\,.
$$
We notice that in order to take the prescribed initial data, $F$ must have the asymptotic behavior$$
F(\xi)/|\xi|^{\gamma} \to C(\theta) \quad \mbox{as } \ |\xi|\to\infty\,.
$$
This means that for every $t>0$
$$
\lim_{|x|\to\infty}F(x)/|x|^{\gamma} =C(\theta).
$$

As an immediate consequence of the existence of these solutions, we derive the optimality of the time exponent 1 in the dependence of formula \eqref{quant.est} for large times.

%%%%%%%%%%%%%%%%%%%%%%%%%%%%%%%%%%%%%%%%%%%%%%%%%%%
\subsection{New weighted estimates}

\noindent {\bf Comparison with known results.} (i) In the Heat Equation  case, $\partial_t u=\Delta u$, the allowed growth for the initial measure is square exponential, and the estimate is quite different. The main qualitative difference  is that  blow-up in finite time may happen. 

(ii) In the Porous Medium Equation, $\partial_t u=\Delta u^m$, $m>1$, the allowed initial growth is roughly $|x|^{2/(m-1)}$ and the estimates are given in \cite{BCP84}, Theorem E, see also \cite{VazPME}. They have a form that is not so different from ours, but again blow-up in finite time may happen. When we compare more closely with our estimate $|x|^{N+2s}$  the exponent $N$ looks strange.

In view of these observations, we look for pointwise bounds with growth at infinity more similar to  \cite{BCP84}. We want to prove that for data with power growth  the solution $u$ keeps pointwise the same growth as the initial datum.

\begin{prop}\label{Weighted.SE}
Let $u$ be a very weak solution of the FHE \eqref{FHE} corresponding to the nonnegative initial datum $0\le u_0\in \mathcal{M}_s^+$. Assume moreover that
\begin{equation}\label{growth.u0}
0\le u_0(x)\le U_0(|x|)\qquad\mbox{for $|x|\ge R_0\ge 1$.}
\end{equation}
for some nondecreasing function $U_0: [0,\infty)\to[0,\infty)$\,.
 Then we have the following bounds:
\begin{equation}\label{Weighted.SE.1}
\frac{ u(t,x) }{\left(1+ U_0(2|x|)\right)}\le  c_2\,\left(t +t^{-N/2s}\right)
\|\mu_0\|_\Phi\,.
\end{equation}
\end{prop}
\noindent\textbf{Remarks. } (i) The radial bound is not strange if we think of the property of symmetrization, \cite{VazVol1}.

(ii) We know that the power function $U_0(x)=|x|^{2s}$ is not admissible as initial datum, but all lower powers  $|x|^\alpha$ with $\alpha<2s$ are admissible.

\noindent {\sl Proof.~}By linearity we can assume $\|\mu_0\|_\Phi=1$. Let us fix $x\in \RR^N$ and split two regions:
\begin{equation}\label{A1.A2}
A_1=\left\{y\in\RR^N\,\big|\, |x-y|\ge |y|/2\right\}\qquad\mbox{and}\qquad A_2=A_1^c=\left\{y\in\RR^N\,\big|\, |x-y|\le |y|/2\right\}
\end{equation}
Since $u$ is a very weak solution, it can be represented with the formula \ref{repres.form} so that
\begin{equation}\label{u.repr}
u(t,x)=\int_{\RR^N}u_0(y)P^t(x-y)\dy
\end{equation}
and thanks to the bounds \eqref{HK-estimates.1} for $P^t$ we have
\begin{equation}\label{HK-estimates.1b1}\begin{split}
u(t,x)&\le c_1\,t \int_{\RR^N}\frac{u_0(y)}{\big(t^{1/s}+|x-y|^2\big)^{(N+2s)/2}}\dy \\
&= c_1\,t \int_{A_1}\frac{u_0(y)}{\big(t^{1/s}+|x-y|^2\big)^{(N+2s)/2}}\dy
+c_1\,t \int_{A_2}\frac{u_0(y)}{\big(t^{1/s}+|x-y|^2\big)^{(N+2s)/2}}\dy = I+II\,.
\end{split}
\end{equation}
Let us estimate the two integrals separately. On one hand, since on $A_1$ we have $|x-y|\ge |y|/2$\,, then we easily get
\begin{equation}\label{HK-estimates.1c}\begin{split}
I &\le c_2\,t \int_{A_1}\frac{u_0(y)}{\big(t^{1/s}+|y|^2\big)^{(N+2s)/2}}\dy
= c_2\,t \int_{A_1}\frac{u_0(y)}{\big(t^{1/s}+|y|^2\big)^{(N+2s)/2}}\left(\frac{1+|y|^2}{1+|y|^2}\right)^{N+2s}\dy\\
&\le c_2\,t \left(t^{-1/2s}+1\right)^{N+2s} \int_{A_1}\frac{u_0(y)}{\big(1+|y|^2\big)^{(N+2s)/2}}\dy \\
&\le c_2\,t \left(t^{-1/2s}+1\right)^{N+2s} \int_{\RR^N}\frac{u_0(y)}{\big(1+|y|^2\big)^{(N+2s)/2}}\dy
\end{split}\end{equation}
where we have used the simple inequality
\begin{equation}\label{simple.ineq}
\frac{1+|y|^2}{t^{1/2s}+|y|^2}=\frac{1 }{t^{1/2s}+|y|^2}+\frac{|y|^2}{t^{1/2s}+|y|^2}\le t^{-1/2s}+1\,.
\end{equation}
On the other hand, on $A_2$\,, we first notice that $|y|\le |x-y|+|x|\le |y|/2 +|x|$\,, from which we deduce that $|y|\le 2|x|$. This latter fact, combined with the fact that $|x-y|\le |y|/2\le |x|$ allows to deduce that $A_2\subseteq B_{|x|}(x)$\,.  Finally, using hypothesis \eqref{growth.u0} we see that\,, if $|x|\ge R_0$ and $y\in A_2$\,, we get
\[
u_0(y)\le U_0(|y|)\le U_0(2|x|)\,,
\]
since $U$ is non-decreasing\,, so that
\begin{equation}\label{HK-estimates.1d}\begin{split}
II &\le c_1\,t \int_{B_{|x|}(x)}\frac{u_0(y)}{\big(t^{1/s}+|x-y|^2\big)^{(N+2s)/2}}\dy\\
&\le c_1\,t U_0(2|x|)\int_{B_{|x|}(x)}\frac{1}{\big(t^{1/s}+|x-y|^2\big)^{(N+2s)/2}}\dy\\
&\le c_1\,t U_0(2|x|)\left(t^{-1/2s}+1\right)^{N+2s} \int_{B_{|x|}(x)}\frac{1}{\big(1+|x-y|^2\big)^{(N+2s)/2}}\dy\\
&\le c_1\,t U_0(2|x|)\left(t^{-1/2s}+1\right)^{N+2s} \int_{\RR^N}\frac{1}{\big(1+|x-y|^2\big)^{(N+2s)/2}}\dy\\
&\le c_2\,t U_0(2|x|)\left(t^{-1/2s}+1\right)^{N+2s}\,.
\end{split}
\end{equation}
where in the last step we have used again inequality \eqref{simple.ineq} in the form
\[
\frac{1+|x-y|^2}{t^{1/2s}+|x-y|^2}=\frac{1 }{t^{1/2s}+|x-y|^2}+\frac{|x-y|^2}{t^{1/2s}+|x-y|^2}\le t^{-1/2s}+1\,.
\]
and the fact that
\[
\int_{\RR^N}\frac{1}{\big(1+|x-y|^2\big)^{(N+2s)/2}}\dy\le C<\infty\,.
\]
Joining then \eqref{HK-estimates.1c} and \eqref{HK-estimates.1d} we conclude the proof.\qed

\section{Regularity properties}\label{sect.regularity}

\subsection{Global Estimates}

We contribute here an interesting result on global estimates which gives pointwise global bounds for solutions to the FHE in terms of the heat kernel, or equivalently  we recall the estimates \eqref{HK-estimates.1}
\begin{equation*}
P^t(x)\asymp  \frac{t}{\big(t^{1/s}+|x|^2\big)^{(N+2s)/2}}\,.
\end{equation*}

\begin{thm}[Global Estimates]\label{thm.GHP}
Let $u$ be a very weak solution of the FHE \eqref{FHE} corresponding to the nonnegative initial datum $0\le \mu_0\in \mathcal{M}_s^+$. Assume moreover that
\begin{equation}\label{growth.u0.L1}
0\le u_0(x)\le \frac{1}{(1+|x|^2)^{\frac{N+2s}{2}}}\qquad\mbox{for $|x|\ge R_0\ge 1$.}
\end{equation}
Then we have the following bounds for all $t>0$ and $x\in \RR^N$
\begin{equation}\label{Glob.Harnack}
k_0 P^t(x)\|u_0\|_{\LL^1(\RR^N)}\le u(t,x)\le k_1 P^t(x)\|u_0\|_{\LL^1(\RR^N)}\,.
\end{equation}
where the constants $k_0,k_1$ depend on $N,s$ and $R_0$\,.
\end{thm}
\noindent {\sl Proof.~}The proof is divided in two steps. By linearity we can assume $\|u_0\|_{\LL^1(\RR^N)}=1$. Since $u$ is a very weak solution, it can be represented with the formula \eqref{repres.form} so that
\begin{equation}\label{u.repr.GHP}
u(t,x)=\int_{\RR^N}u_0(y)P^t(x-y)\dy\,.
\end{equation}
By a rescaling argument, cf. the end of Step 1 for more details, it is sufficient to prove our bounds at time $t=1$.

\noindent$\bullet~$\textsc{Step 1. }\textit{Upper bounds. }Let us fix $x\in \RR^N$ and split three regions:
\begin{equation}\label{A1.A2.A3}
A_1=B_{|x|/3}(0)\,,\qquad A_2= B_{|x|/3}^c(0)\setminus B_{|x|/10}(x)
\qquad\mbox{and}\qquad A_3=B_{|x|/10}(x)\,.
\end{equation}
We will estimate the right-hand side of \eqref{u.repr.GHP} on the regions $A_1,A_2,A_3$ separately.

On $A_1$ we have that $|x-y|\ge |x|-|y|\ge 2|x|/3$, and thanks to the bounds \eqref{HK-estimates.1} for $P^t$ we have
\begin{equation}\label{A1.up.GHP}
\int_{A_1}u_0(y)P^1(x-y)\dy
\le c_1 \int_{A_1}\frac{u_0(y)}{(1+|x-y|^2)^{\frac{N+2s}{2}}}\dy
\le c_1\frac{\int_{A_1}u_0(y)\dy}{(1+\frac{4}{9}|x|^2)^{\frac{N+2s}{2}}}
\le \frac{c_2}{(1+ |x|^2)^{\frac{N+2s}{2}}}\,,
\end{equation}
we have used that $\|u_0\|_{\LL^1(\RR^N)}=1$\,.

On $A_2$ we have that $|x-y|\ge |x|/10$, and thanks to the bounds \eqref{HK-estimates.1} for $P^t$ we have
\begin{equation}\label{A2.up.GHP}
\int_{A_2}u_0(y)P^1(x-y)\dy
\le c_1 \int_{A_2}\frac{u_0(y)}{(1+|x-y|^2)^{\frac{N+2s}{2}}}\dy
\le c_1\frac{\int_{A_2}u_0(y)\dy}{(1+\frac{1}{10}|x|^2)^{\frac{N+2s}{2}}}
\le \frac{c_3}{(1+ |x|^2)^{\frac{N+2s}{2}}}
\end{equation}
we have used that $\|u_0\|_{\LL^1(\RR^N)}=1$\,.

On $A_3$ we have that $|y|\ge |x|-|x-y|\ge 9|x|/10$, and we use that $|x|\ge R_0$\,, so that we can use the assumption on the initial datum $u_0(y)\le (1+ |x|^2)^{-\frac{N+2s}{2}}$. We notice that this can not be avoided in view of the possibility to build counter examples in the spirit of Section \ref{sec.boundedness}.  Therefore for all $y\in A_3$ we have that
\begin{equation}\label{A3.0.up.GHP}
u_0(y)\le  \frac{1}{(1+ |y|^2)^{\frac{N+2s}{2}}}\le   \frac{1}{(1+\frac{9}{10}|x|^2)^{\frac{N+2s}{2}}}\le   \frac{c'_4}{(1+ |x|^2)^{\frac{N+2s}{2}}}
\end{equation}
Thanks to the bounds \eqref{HK-estimates.1} for $P^t$ we have
\begin{equation}\label{A3.up.GHP}\begin{split}
\int_{A_3}  u_0(y)P^1(x-y)\dy
&\le c_1 \int_{A_3}\frac{u_0(y)}{(1+|x-y|^2)^{\frac{N+2s}{2}}}\dy\\
&\le \frac{c'_4}{(1+ |x|^2)^{\frac{N+2s}{2}}} \int_{A_3}\frac{1}{(1+|x-y|^2)^{\frac{N+2s}{2}}}\dy
\le \frac{c_4}{(1+ |x|^2)^{\frac{N+2s}{2}}}
\end{split}\end{equation}
we have used that
\[
\int_{A_3}\frac{1}{(1+|x-y|^2)^{\frac{N+2s}{2}}}\dy\le \int_{\RR^N}\frac{1}{(1+|z|^2)^{\frac{N+2s}{2}}}\rd z\le c''_4<+\infty\,.
\]
Finally, joining \eqref{A1.up.GHP},  \eqref{A2.up.GHP} and \eqref{A3.up.GHP} and undoing the renormalization $\|u_0\|_{\LL^1(\RR^N)}=1$\,, we can estimate the right-hand side of \eqref{u.repr.GHP} as follows
\begin{equation}\label{step1.up.GHP.t=1}
u(1,x)=\int_{\RR^N}u_0(y)P^1(x-y)\dy \le \frac{c_2+c_3+c_4}{(1+ |x|^2)^{\frac{N+2s}{2}}}\|u_0\|_{\LL^1(\RR^N)}
:= c_5\|u_0\|_{\LL^1(\RR^N)} P^1(x)
\end{equation}
By standard rescaling, applying the above inequality to $u_\lambda(1,x)= u(\lambda , \lambda^{-1/2s}x)$)  then letting $\lambda=t$\,, and since we know that $P^t(x)=P^1\left(xt^{-1/2s}\right)t^{-N/2s}\asymp \big(1+|x|^2t^{-1/s}\big)^{-(N+2s)/2} t^{-N/2s}$ we obtain
\begin{equation}\label{step1.up.GHP}
u(t,x)=\int_{\RR^N}u_0(y)P^t(x-y)\dy \le \frac{c_5 \|u_0\|_{\LL^1(\RR^N)}}{(1+ \frac{|x|^2}{t^{1/s}})^{\frac{N+2s}{2}}}\frac{1}{t^{N/2s}}
:= k_1\|u_0\|_{\LL^1(\RR^N)} P^t(x)\,,
\end{equation}
which proves the upper bound in \eqref{Glob.Harnack}.

\noindent$\bullet~$\textsc{Step 2. }\textit{Lower bounds. }Let $M=\int_{\RR^N}u_0\dy>0$. Since $\int_{B_R(0)}u_0\dy\to M$ as $R\to \infty$\,, then there always exist $R_M>0$ such that $\int_{B_{R_M(0)}}u_0\dy\ge M/2$. Notice that $|x-y|\le |x|+|y|\le |x|+R_M$, so that the bounds on the heat kernel \eqref{HK-estimates.1}  and the representation formula \eqref{u.repr.GHP} imply
\begin{equation}\label{u.repr.low.GHP}\begin{split}
u(1,x)&=\int_{\RR^N}u_0(y)P^1(x-y)\dy \ge c_1 \int_{B_{R_M(0)}}\frac{u_0(y)}{(1+|x-y|^2)^{\frac{N+2s}{2}}}\dy\\
&\ge c_2\frac{\int_{B_{R_M(0)}}u_0(y)\dy}{(1+ (|x|+R_M)^2)^{\frac{N+2s}{2}}}
 \ge \frac{M}{2}\frac{c_2}{(1+ (|x|+R_M)^2)^{\frac{N+2s}{2}}}\\
 &:= \frac{k_0(R_M)\, M}{(1+ |x|^2)^{\frac{N+2s}{2}}}
 =k_0(R_M)\|u_0\|_{\LL^1(\RR^N)}P^1(x)\,.
\end{split}\end{equation}
The same scaling argument as at the end of Step 1, finally gives
\begin{equation}\label{step2.low.GHP}
u(t,x)=\int_{\RR^N}u_0(y)P^t(x-y)\dy \ge \frac{k_0(R_M) \|u_0\|_{\LL^1(\RR^N)}}{(1+ \frac{|x|^2}{t^{1/s}})^{\frac{N+2s}{2}}}\frac{1}{t^{N/2s}}
= k_0(R_M)\|u_0\|_{\LL^1(\RR^N)} P^t(x)\,,
\end{equation}
which proves the lower bound in \eqref{Glob.Harnack}.\qed

\subsection{Harnack inequalities and H\"older regularity}

In this section, we deal with Harnack inequalities and H\"older regularity for the solutions of the FHE. Actually, there has been lately a substantial literature dealing with nonlocal parabolic equations (see for instance \cite{CCV,FK,KS} and references therein). As far as the FHE is concerned, one can quote the paper  \cite{BdP} where smoothing effects and decay properties are shown. In the previous papers, the authors deal also with more general kernels than the one of the fractional Laplacian.

As a corollary of the global estimates of the previous section, we obtain a quite precise local Harnack inequality.
We recall that the first Harnack inequality for the standard heat equation is due to Pini \cite{Pi54} and Hadamard \cite{Ha54} and reads
\[
u(t_1,x)\le u(t_2,y)\left(\frac{t_2}{t_1}\right)^{N/2}\ee^{\frac{|x-y|^2}{4(t_2-t_1)}}\qquad\mbox{for any }x,y\in \RR^N t_2>t_1>0\,.
\]
The above inequality was then generalized by Moser \cite{Mo64} to solutions to more general linear parabolic operators on the cylinder $Q=(0,T)\times B_{R}$
\[
\sup_{Q_-}u\le \mathcal{H}\inf_{Q_+}u
\]
for any $Q_-=[t_1,t_2]\times B_r(x_0) $\,, $Q_-= [t_3,t_4]\times B_r(x_0)$\,, such that $t_2<t_3$ and $0<r<R$.  The constant $\mathcal{H}$ depends on $N, t_1,\dots, t_4, R, r$\,. Notice that this is called forward Harnack inequality and $\delta=t_3-t_2>0$ is a delay. One may think that the delay can be neglected, but it is can not for the classical heat equation. For a counterexample we need only take displaced copies of the fundamental solution, the Gaussian, and compare supremum and infimum on the same ball and at the same time is impossible, hence the delay is really needed.

We show here that a stronger Harnack inequality holds for the FHE. Indeed the delay can be neglected as we will show in the following Theorem.

\begin{thm}[\bf Local Harnack inequalities of Forward/Backward/Elliptic type]
Under the assumptions of Theorem $\ref{thm.GHP}$\,, we have the following Harnack inequality for all $x,y\in \RR^N$ and all $t,\tau>0$
\begin{equation}\label{Local.Harnack.1}
u(t,x)\le \mathcal{C}_1 \frac{t}{\tau}
\left[1+\frac{|t-\tau|^{\frac{1}{s}}+\big||y|^2-|x|^2\big|}{t^{\frac{1}{s}}+|x|^2}\right]^{\frac{N+2s}{2}} u(\tau,y)
%\le \mathcal{C}_1 \left(\frac{\tau}{t}\right)^{\frac{N}{2s}}
%\left[\frac{1+|y|^2\tau^{-\frac{1}{s}}}{1+|x|^2t^{-\frac{1}{s}}}\right]^{\frac{N+2s}{2}} u(\tau,y)
\end{equation}
Moreover, the backward/forward/elliptic Harnack inequality holds, namely for all $0<t,\tau\le T$
\begin{equation}\label{Local.Harnack.2}
\sup_{x\in B_R(x_0)}u(t,x)\le \mathcal{C}_2 \inf_{y\in B_R(x_0)} u(\tau,y)\,.
\end{equation}
The constants $\mathcal{C}_1\sim k_1/k_0$, where $k_0,k_1$ are as in Theorem $\ref{thm.GHP}$\,. On the other hand, in general $\mathcal{C}_2=\mathcal{C}_2(N,s,\tau,t,T,R)$. It has an explicit bound given in formula \eqref{Harnack.loc.1}. When  $R$ is ``small'' compared to $|x_0|$, then $\mathcal{C}_2$ simplifies to
\begin{equation}\label{Local.Harnack.2b}
\mathcal{C}_2=\mathcal{C}_1 \frac{t}{\tau}\left[19+\left(\frac{T}{t}\right)^{\frac{1}{s}} \right]^{\frac{N+2s}{2}}\qquad\mbox{for all $R<|x_0|/2$}\,.
\end{equation}
\end{thm}

\noindent\textbf{Remark.}  This stronger result may seem surprising, but it is known in nonlinear diffusions, where it is typical to have ``reversed in time'' Harnack inequalities, especially in the regime of Fast Diffusion Equations, cf. \cite{BV-ADV, DGVbook}.

\noindent {\sl Proof.~}Formula \eqref{Local.Harnack.1} follows by the global bounds \eqref{Glob.Harnack}, simply by recalling the heat kernel estimates $P^t(x) \asymp t \big(t^{1/s}+|x|^2\big)^{-(N+2s)/2}$. The constant $\mathcal{C}_1\sim k_1/k_0$\,, where $k_0,k_1$ are as in Theorem $\ref{thm.GHP}$\,. Indeed, we can assume without loss of generality that $\|u_0\|_{\LL^1(\RR^N)}=1$. Recall that   so that
\begin{equation}\label{Harnack.Pt.0}
P^t(x) \le \frac{c_1\, t}{\big(t^{\frac{1}{s}}+|x|^2\big)^{\frac{N+2s}{2}}}\qquad\mbox{and that }\qquad  P^{\tau}(y) \ge \frac{c_0\, \tau}{\big(\tau^{\frac{1}{s}}+|y|^2\big)^{\frac{N+2s}{2}}}
\end{equation}
so that
\begin{equation}\label{Harnack.Pt.1}
P^t(x)\le \frac{c_1\, t}{c_0\tau}\frac{\big(\tau^{\frac{1}{s}}+|y|^2\big)^{\frac{N+2s}{2}}}{\big(t^{\frac{1}{s}}+|x|^2\big)^{\frac{N+2s}{2}}}P^{\tau}(y)
\end{equation}
from which \eqref{Local.Harnack.1} follows simply by using \eqref{Glob.Harnack}\,.

We next notice that if $x,y\in B_R(x_0)$ and $0<t,\tau \le T$\,, then
\begin{equation}\label{Harnack.loc.0}
\mathcal{C}_1 \frac{t}{\tau}\left[1+\frac{|t-\tau|^{\frac{1}{s}}+\big||y|^2-|x|^2\big|}{t^{\frac{1}{s}}+|x|^2}\right]^{\frac{N+2s}{2}}
\le \mathcal{C}_1 \frac{t}{\tau}\left[1+\frac{T^{\frac{1}{s}}+2(|x_0|+R)^2}{t^{\frac{1}{s}}+(|x_0|-R)^2}\right]^{\frac{N+2s}{2}}
\end{equation}
This provides a finite expression for $\mathcal{C}_2$. Moreover, under the condition $R<|x_0|/2$ we can estimate better the last term, and eliminated the dependence on $R$ and $x_0$\,, as follows
\begin{equation}\label{Harnack.loc.0b}\begin{split}
\left[1+\frac{|t-\tau|^{\frac{1}{s}}+\big||y|^2-|x|^2\big|}{t^{\frac{1}{s}}+|x|^2}\right]^{\frac{N+2s}{2}}
&\le \left[1+\left(\frac{T}{t}\right)^{\frac{1}{s}}+2\frac{(|x_0|+R)^2}{(|x_0|-R)^2}\right]^{\frac{N+2s}{2}}\\
&\le \left[1+\left(\frac{T}{t}\right)^{\frac{1}{s}}+18\frac{|x_0|^2}{|x_0|^2}\right]^{\frac{N+2s}{2}}
=\left[19+\left(\frac{T}{t}\right)^{\frac{1}{s}} \right]^{\frac{N+2s}{2}}
\end{split}
\end{equation}
which justifies the bound \eqref{Local.Harnack.2b}\,.
Once \eqref{Local.Harnack.1} is established and its constant estimated as above, we just take supremum in $x\in B_R(x_0)$ and the infimum in $y\in B_R(x_0)$ to obtain
\begin{equation}\label{Harnack.loc.1}
\sup_{x\in B_R(x_0)}u(t,x)\le  \mathcal{C}_1 \frac{t}{\tau}
\left[1+\frac{T^{\frac{1}{s}}+2(|x_0|+R)^2}{t^{\frac{1}{s}}+(|x_0|-R)^2}\right]^{\frac{N+2s}{2}} \inf_{y\in B_R(x_0)} u(\tau,y)
\end{equation}
which proves \eqref{Local.Harnack.2}\,.\qed

\noindent {\bf Related literature.} In the seminal paper \cite{CCV}, the authors investigate a linear nonlocal parabolic equation with a general kernel bounded from above and from below by the one of the fractional Laplacian and prove H\"older regularity of the solutions following the approach of De Giorgi. Using (weak) Harnack inequalities, the H\"older continuity of the solutions has been proved in \cite{FK, KS} for instance.

For the case of fully nonlinear (i.e. in non-divergence form) nonlocal equations we refer to \cite{CLD2,CLD3} for H\"older continuity results. In the case of (possibly degenerate/singular)nonlinear diffusion, it is a difficult problem to prove full regularity. In \cite{AthC}, the authors prove H\"older regularity for the nonlinear fractional diffusion equation in the non-degenerate case. In the paper \cite{VPQR}, a thorough study of the regularity of the solutions of the fractional filtration equation (and hence fractional PME) is investigated. In the linear case (with a right hand side), the authors prove that the solutions are classical and smooth.

\subsection{Analiticity and Gevrey regularity}\label{sec.gevrey}
As already mentioned, the fractional heat kernel has a smooth Fourier symbol, $\ee^{-|\xi|^{2s}}$. The sell-known relation between the spatial decay of the Fourier transform and regularity of the original function, immediately shows that $P_t$ is also a smooth function, in general $C^\infty$. The regularity of the solutions is therefore the same, in virtue of the representation formula $u=P^t\ast u_0$, and standard properties of convolution.

The next question that arises, is about finer regularity properties for solutions of the FHE, like analyticity. In the classical case, $s=1$\,, solutions to the HE are known to be analytic, since the Gaussian kernel is analytic, and this can be checked directly since it is explicit. In the case of the FHE the panorama is a bit different, the kernel is not explicit. To explain the finer regularity properties of solution to the FHE, it is convenient to recall some concepts.

We recall that a function $f\in \LL^1(\RR^N)$ belongs to the Gevrey class $G^{1/2s}$ if and only if there exists $\varepsilon>0$ such that
\begin{equation}\label{gevrey.def}
\ee^{\varepsilon (1+|\xi|^2)^s}\hat{f}(\xi)\in \LL^2(\RR^N)\,.
\end{equation}
Notice that $G^s\subset \bigcap_{\sigma\ge 0}H^\sigma(\RR^N)\subset C^\infty(\RR^N)$. The class of Gevrey function can also be characterized via Cauchy inequalities, as the class of functions such that the $k^{th}$ derivatives do not grow more than $(k!)^{1/2s}$, of course if $1/2s\le 1$ this functions are analytic, since the $k^{th}$ derivatives do not grow more than $k!$, indeed when $1/2s< 1$ this class is contained in the class of the so-called ultra-analytic functions\,. On the other hand, when $1/2s> 1$ there may be function in the Gevrey class which are not analytic, in particular nontrivial compactly supported functions (clearly not analytic), cf. for instance \cite{rodino}.

The optimal regularity for solutions to the FHE $u_t+(-\Delta)^s u=0$\,, then reads: for any $u_0\in \LL^1$\,, the corresponding solution $u(t,\cdot)$ satisfies for every $t>0$
\begin{equation}\label{Analiticity.Gevrey.1}
u(t, \cdot)=P^t\ast u_0\in G^{\frac{1}{2s}}(\RR^N)
%\left\{\begin{array}{lll}
%\subset \mathcal{A} & \mbox{if $\frac{1}{2}<s<1$ solutions are ultra-analytic}\\
%=\mathcal{A}& \mbox{if $s=\frac{1}{2}$ solutions are analytic (explicit kernel)}\\
%=G^{1+\frac{1-2s}{2s}}& \mbox{if $0<s<\frac{1}{2}$ solutions are not analytic, only Gevrey}\\
%\end{array}\right.
\end{equation}
and we recall the panorama
\begin{equation}\label{Analiticity.Gevrey.2}
G^{\frac{1}{2s}}(\RR^N)=
\left\{\begin{array}{lll}
\subset \mathcal{A}(\RR^N),& \mbox{if $\frac{1}{2}<s<1$,}& \mbox{solutions are ultra-analytic}\\
=\mathcal{A}(\RR^N),& \mbox{if $s=\frac{1}{2}$,} & \mbox{solutions are analytic (explicit kernel)}\\
=G^{1+\frac{1-2s}{2s}}(\RR^N),& \mbox{if $0<s<\frac{1}{2}$,} & \mbox{solutions are not analytic, only Gevrey.}\\
\end{array}\right.
\end{equation}

%%%%%%%%%%%%%%%%%%%%%%%%%%%%%%%%%%%%%%%

\section{Theory for very weak solutions with two signs}\label{sec.signed}

We have restricted our data and solutions to be nonnegative in order to obtain an optimal theory. However, a large part of the results are true for signed solutions. First of all, by splitting the initial measure $\mu_0$ into its positive and negative parts and using the linearity and the representation formula \eqref{repres.form} separately on each part, we can obtain a formal solution for data $\mu_0\in \mathcal{M}_s(\RR^N)$ defined as class of locally finite Radon measures
\begin{equation}
\int_{\ren} (1+|x|)^{-(N+2s)}\,d|\mu|(x)<\infty\,.
\end{equation}
It is immediate that for every $\mu_0\in \mathcal{M}_s(\RR^N)$ we obtain a signed very weak solution in the sense of Definition \ref{def.vweak} without the restriction of positivity. The initial data are taken in the sense of measures as in \eqref{init.trace.def}.

\noindent$\bullet~$ An important and nontrivial result is that uniqueness holds also for signed solutions with this definition of very weak solution. This is the merit of the proof of Theorem \ref{thm.uniq}.\\
 It is clear from the proof that only minor changes are needed to prove Theorem \ref{thm.uniq} for sign-changing solutions. Indeed, notice that in the original proof of the theorem the function $u$ which is by definition the difference of two very weak solutions has no definite sign. Therefore, most of the proof remains the same, except that the property $|u|=u_1+u_2$ has to be replaced by the trivial bound $|u| \leq |u_1|+|u_2|$ in the computation of integral $J_{1,n}$. Finally, one needs to check that Lemma \ref{lem.weighted.L1} holds for sign-changing solutions, but this is a direct consequence of the linearity of the equation.

We summarize the above discussion in the following theorem on existence and uniqueness for signed solution with data in $\mathcal{M}_s$.

\begin{thm}\label{thm.exist.uniq.signed} For every initial data $\mu_0 \in \mathcal{M}_s$.
there exists a unique very weak solution to the FHE \eqref{FHE} in the sense of above definition that is given by the representation formula \eqref{repres.form}.
\end{thm}

The classical heat equation theory tells us that there may be oscillating solutions with a higher growth at infinity than the growth of nonnegative solutions. For such new classes of highly growing signed solutions, events like the Tychonov non-uniqueness theorem may happen, cf. \cite{Ty}.

\medskip

We cannot prove the theorem of existence of initial traces for signed very weak solutions. This is a common difficulty with signed solutions. This means that for signed solutions the theory is not optimal.

\section{Comments and extensions}

\noindent$\bullet~$  {\bf The equation with right-hand side.} The fact that the evolution semigroup associated to the FHE is contractive can be seen as a consequence of the fact that the fractional Laplacian operator is $m$-accretive in all spaces $L^p(\ren)$, $1\le p\le \infty$, the same as the classical heat equation.
The Crandall-Liggett theorem \cite{CL71} will then produce a solution of
%Why not for $f\in L^1(0,T:L^1_\Phi(\ren))$. We have the Duhamel formula
\begin{equation}\label{FHEf}
\partial_t u + (-\Delta)^s u=f(x,t), \quad 0<s<1\,.
\end{equation}
for data $f\in L^1(0,T:L^1(\ren))$.

\medskip

%\noindent$\bullet~$ The question of natural growth that is proved by Barrios et al for the good solutions that they call strong. Comment on viscosity solutions maybe not needed.
%
%\medskip

%\noindent$\bullet~$ Other examples of explicit solutions

\noindent$\bullet~$ There are many examples of diffusion equations of nonlocal type, that are investigated at this moment, probably some of our results will apply.

- Parabolic equations with fractional derivative in time by Allen, Caffarelli and Vasseur \cite{ACV1}. %,ACV2}

- Extension property for a parabolic operator Stinga and Torrea  \cite{ST}.

%\newpage

\section{Appendix}
We recall here Lemma 2.1 of \cite{BV2012} since we use it several times in the paper.
\begin{lem}\label{Lem.phi}
Let $\varphi\in C^2(\RR^N)$ be a positive real function that is radially symmetric and decreasing in $|x|\ge 1$. Assume also that $\varphi(x)\le |x|^{-\alpha}$ and that $|D^2\varphi(x)| \le c_0 |x|^{-\alpha-2}$\,, for some positive constant $\alpha$ and for $|x|$ large enough. Then,  for all $|x|\ge |x_0|>>1$ we have
\begin{equation}\label{Delta.s.phi}
|(-\Delta)^s\varphi(x)|\le
\left\{\begin{array}{lll}
\dfrac{c_1}{|x|^{\alpha+2s}}\,,   & \mbox{if $\alpha<N$}\,,\\[5mm]
\dfrac{c_2\log|x|}{|x|^{N+2s}}\,,   & \mbox{if $\alpha=N$}\,,\\[5mm]
\dfrac{c_3}{|x|^{N+2s}}\,,   & \mbox{if $\alpha>N$}\,,\\[5mm]
\end{array}\right.
\end{equation}
with positive constants $c_1,c_2,c_3>0$ that depend only on $\alpha,s,N$ and $\|\varphi\|_{C^2(\RR^N)}$. For  $\alpha>d$ the reverse estimate holds  from below  if $\varphi\ge0$: \ $|(-\Delta)^s\varphi(x)|\ge c_4 |x|^{-(N+2s)}$  for all $|x|\ge |x_0|>>1$\,.
\end{lem}

%%%%%%%%%%%%%%%%%%%%%%%%%%%%%%%%%%%%%%%%%%%%%%%%%%%%%%%%%%%%%%%%%%%%%%%
\
%\vfill

\noindent {\large \sc Acknowledgment}

\noindent  First and third author funded by Project MTM2014-52240-P  (Spain).  The second author would like to thank the hospitality of the Departamento de Matem\'{a}ticas of Universidad Aut\'{o}noma de Madrid where this work started. MB and JLV would like to thank the hospitality of the Departments of Mathematics of  the Johns Hopkins University and of the University of Texas at Austin for their kind hospitality.

%%%%%%%%%%%%%%%%%%%%%%%%%%%%%%%%%%%%%%%%%%%%%%%%%%%%%%%%%%%%%%%%%%%%%
\

%\newpage

%%%%%%%%%%%%%%%%%%%%%%%%%%%%%%%%%%%%%%%%%%%%%%%%%%%%%%%%%%%%%%%%%%%%%

\bibliographystyle{amsplain} 

\

\noindent {\bf Keywords.} Fractional Heat Equation, Existence and uniqueness theory, Optimal data, Pointwise Estimates, Harnack inequalities.

\medskip

\noindent {\sc Mathematics Subject Classification}.   35A01, 35A02, 35K99
 
\noindent e-mail addresses:\texttt{~matteo.bonforte@uam.es,~sire@math.jhu.edu,~juanluis.vazquez@uam.es}

\end{document}